\documentclass[reqno]{amsart}

\usepackage[latin1]{inputenc}
\usepackage{subfigure}
\usepackage{pstricks}
\usepackage{fp}
\usepackage{multido}
\usepackage{amsmath}

\newtheorem{theo}{Theorem}
\newtheorem{prop}{Proposition}

\newtheorem{lemm}[theo]{Lemma}

\newtheorem{defi}[theo]{Definition}

\newtheorem{rema}[theo]{Remark}

\newcommand{\pref}[1]{(\ref{#1})}

\def\cf{{\it cf. }}
\def\ie{{\it i.e. }}

\def\si{\sigma}

\def\bx{\bar x}

\def\ba{\bar a}

\newcommand{\flecheD}{\psline{->}(0,0)(.75,0)\psline(0,0)(1,0)}
\newcommand{\flecheG}{\psline{->}(0,0)(-.75,0)\psline(0,0)(-1,0)}
\newcommand{\flecheH}{\psline{->}(0,0)(0,.75)\psline(0,0)(0,1)}
\newcommand{\flecheB}{\psline{->}(0,0)(0,-.75)\psline(0,0)(0,-1)}

\newcounter{ISL}
\newcounter{ISH}
\newcommand{\IceGrid}[2]{%
\setcounter{ISL}{#1}\addtocounter{ISL}{-1}%
\setcounter{ISH}{#2}\addtocounter{ISH}{-1}%
\multido{\i=0+1}{#1}{\psline(\i,0)(\i,\theISH)}%
\multido{\i=0+1}{#2}{\psline(0,\i)(\theISL,\i)}%
}

\newcommand{\colD}[1]{%
\multido{\i=0+1}{#1}{\rput(0,\i){\flecheD}}%
}

\newcommand{\colG}[1]{%
\multido{\i=0+1}{#1}{\rput(1,\i){\flecheG}}%
}

\newcommand{\linH}[1]{%
\multido{\i=0+1}{#1}{\rput(\i,0){\flecheH}}%
}

\newcommand{\linB}[1]{%
\multido{\i=0+1}{#1}{\rput(\i,1){\flecheB}}%
}

\newcommand{\col}[1]{%
\multido{\i=0+1}{#1}{\rput(0,\i){\psline(0,0)(1,0)}}%
}

\newcommand{\lin}[1]{%
\multido{\i=0+1}{#1}{\rput(\i,0){\psline(0,0)(0,1)}}%
}

\newcommand{\halfCC}[5]{%
\degrees[360]%
\multido{\n=#3+#4}{#5}{\rput(#1,#2){\psarc(0,0){\n}{-90}{90}}}%
}

\newcommand{\tqCC}[5]{%
\degrees[360]%
\multido{\n=#3+#4}{#5}{\rput(#1,#2){\psarc(0,0){\n}{-90}{180}}}%
}

\newcommand{\IceSquare}[1]{%
\rput(1,1){\IceGrid{#1}{#1}}%
\rput(0,1){\colD{#1}}%
\rput(1,0){\linB{#1}}%
\rput(#1,1){\colG{#1}}%
\rput(1,#1){\linH{#1}}%
}

\newcounter{ISdoublesize}
\newcommand{\HTIceEven}[1]{%
\setcounter{ISdoublesize}{#1}\addtocounter{ISdoublesize}{#1}%
\rput(1,1){\IceGrid{#1}{\theISdoublesize}}%
\rput(0,1){\colD{\theISdoublesize}}%
\rput(1,0){\linB{#1}}%
\rput(1,\theISdoublesize){\linH{#1}}%
\rput(#1,#1){\rput(0,.5){\halfCC{0}{0}{.5}{1.0}{#1}}}%
\rput(#1,#1){\psline[linestyle=dotted](.25,.5)(.75,.5)}%
}

\newcommand{\HTIceOdd}[1]{%
\setcounter{ISdoublesize}{#1}\addtocounter{ISdoublesize}{#1}%
\addtocounter{ISdoublesize}{1}%
\rput(1,1){\IceGrid{#1}{\theISdoublesize}}%
\rput(0,1){\colD{\theISdoublesize}}%
\rput(1,0){\linB{#1}}%
\rput(1,\theISdoublesize){\linH{#1}}%
\rput(#1,0){\psline(1,1)(1,#1)\rput(1,1){\flecheB}}%
\psarc(#1,#1){1}{0}{90}%
\rput(#1,#1){\halfCC{1}{1}{1}{1}{#1}}%
\rput(#1,1){\col{#1}}%
\rput(#1,#1){\rput(0,2){\col{#1}}}%
\rput(#1,#1){\psline[linestyle=dotted](.25,.25)(1.25,1.25)}
}

\newcounter{QTsize}
\newcommand{\QTIce}[1]{%
\setcounter{QTsize}{#1}\addtocounter{QTsize}{-1}%
\rput(1,1){\IceGrid{\theQTsize}{\theQTsize}}%
\rput(0,1){\colD{#1}}%
\rput(1,0){\linB{#1}}%
\rput(\theQTsize,1){\col{\theQTsize}}%
\rput(1,\theQTsize){\lin{\theQTsize}}%
\rput(\theQTsize,\theQTsize){\psarc(0,0){1}{0}{90}}%
\psline(1,#1)(\theQTsize,#1)\psline(#1,1)(#1,\theQTsize)%
\rput(#1,#1){\tqCC{0}{0}{1}{1}{\theQTsize}}%
\rput(#1,#1){\SpecialCoor\multido{\i=1+1}{\theQTsize}{\rput(\i;45){\psdots[dotstyle=*](0,0)}}}%
\rput(\theQTsize,\theQTsize){\SpecialCoor\psdots[dotstyle=*](1;45)\psline[linestyle=dotted](.5;45)(1.5;45)}%
}

\newcommand{\qQTIce}[1]{%
\setcounter{QTsize}{#1}\addtocounter{QTsize}{-1}%
\rput(1,1){\IceGrid{\theQTsize}{\theQTsize}}%
\rput(0,1){\colD{#1}}%
\rput(1,0){\linB{#1}}%
\rput(\theQTsize,1){\col{\theQTsize}}%
\rput(1,\theQTsize){\lin{\theQTsize}}%
\rput(\theQTsize,\theQTsize){\psarc(0,0){1}{0}{90}}%
\psline(1,#1)(\theQTsize,#1)\psline(#1,1)(#1,\theQTsize)%
\rput(#1,#1){\tqCC{0}{0}{1}{1}{\theQTsize}}%
\rput(#1,#1){\SpecialCoor\multido{\i=1+1}{\theQTsize}{\rput(\i;45){\psdots[dotstyle=*](0,0)}}}%
\rput(\theQTsize,\theQTsize){\SpecialCoor\psline[linestyle=dotted](.5;45)(1.5;45)}%
}

\newcommand{\Hcrossing}{%
\psbezier(0,0)(.4,0)(.6,1)(1,1)%
\psbezier(0,1)(.4,1)(.6,0)(1,0)%
}

\newcommand{\convcorner}{%
\begin{pspicture}(0.2,0.2)\psset{linewidth=.4pt,arrowsize=2.5pt}
\psline{->}(0,.2)(.2,.2)
\psline{->}(.2,0)(.2,.2)
\end{pspicture}%
}

\newcommand{\divcorner}{%
\begin{pspicture}(.2,.2)\psset{linewidth=.4pt,arrowsize=2.5pt}
\psline{->}(.2,.2)(.2,0)
\psline{->}(.2,.2)(0,.2)
\end{pspicture}%
}

\newcommand{\upleft}{%
\begin{pspicture}(.2,.2)\psset{linewidth=.4pt,arrowsize=2.5pt}
\psline(.2,0)(.2,.15)\psarc(.15,.15){.05}{0}{90}\psline{->}(.15,.2)(0,.2)%
\end{pspicture}%
}

\newcommand{\downright}{%
\begin{pspicture}(.2,.2)\psset{linewidth=.4pt,arrowsize=2.5pt}
\psline(0,.2)(.15,.2)\psarc(.15,.15){.05}{0}{90}\psline{->}(.2,.15)(.2,0)%
\end{pspicture}
}

\newcommand{\upc}{%
\begin{pspicture}(.2,.2)\psset{linewidth=.4pt,arrowsize=2.5pt}%
\psline(0,0)(.1,0)\psarc(.1,.1){.1}{270}{90}\psline{->}(.1,.2)(0,.2)%
\end{pspicture}%
}

\newcommand{\downc}{%
\begin{pspicture}(.2,.2)\psset{linewidth=.4pt,arrowsize=2.5pt}%
\psline(0,.2)(.1,.2)\psarc(.1,.1){.1}{270}{90}\psline{->}(.1,0)(0,0)%
\end{pspicture}%
}

\newcommand{\ZqConv}{%
Z_{\textsc{QT}}^{\convcorner}%
}

\newcommand{\ZqDiv}{%
Z_{\textsc{QT}}^{\divcorner}%
}

\newcommand{\ZqUpleft}{%
Z_{\textsc{QT}}^{\upleft}%
}

\newcommand{\ZqDownright}{%
Z_{\textsc{QT}}^{\downright}%
}

\newcommand{\ZhUp}{%
Z_{\textsc{HT}}^{\upc}%
}

\newcommand{\ZhDown}{%
Z_{\textsc{HT}}^{\downc}%
}

\newcommand{\ZhUpleft}{%
Z_{\textsc{HT}}^{\upleft}%
}

\newcommand{\ZhDownright}{%
Z_{\textsc{HT}}^{\downright}%
}

\author{Jean-Christophe~Aval}
\address[J.-C. Aval]
{LaBRI, Universit\'e Bordeaux 1, CNRS\\ 
351 cours de la Lib\'eration,
 33405 Talence cedex, FRANCE}
\email{aval@labri.fr}
\urladdr{http://www.labri.fr/perso/aval}

\author{Philippe Duchon} 
\address[P. Duchon]
{LaBRI, Universit\'e Bordeaux 1, CNRS\\ 
351 cours de la Lib\'eration, 
33405 Talence cedex, FRANCE}
\email{duchon@labri.fr}
\urladdr{http://www.labri.fr/perso/duchon}

\thanks{Both authors are supported by the ANR project MARS (BLAN06-2$\_$0193)}

\title[Quarter-turn Symmetric Alternating Sign Matrices]{Enumeration of alternating sign matrices of even size (quasi-)invariant under a quarter-turn rotation}

\date{\today}

\begin{document}
\maketitle

\begin{abstract}
The aim of this work is to enumerate alternating sign matrices (ASM) that are quasi-invariant under a quarter-turn.
The enumeration formula (conjectured by Duchon) involves, as a product of three terms, the number of unrestricted ASM's and the number of half-turn symmetric ASM's.
\end{abstract}

\section{Introduction}\label{sec:intro}

An {\em alternating sign matrix} is a square matrix with entries in $\{-1,0,1\}$ and such that in any row and column: the non-zero entries alternate in sign, and their sum is equal to $1$. Their enumeration formula was conjectured by Mills, Robbins and Rumsey \cite{MRR}, and proved years later by Zeilberger \cite{zeil}, and almost simultaneously by Kuperberg \cite{kup1}.

Kuperberg's proof is based on the study of the partition function of a
square ice model whose states are in bijection with ASM's. Kuperberg was able
to get an explicit formula for the partition function for some special
values of the spectral parameter. To do this, he used a Yang-baxter
formula and recursive relations obtained by Korepin \cite{korepin} for
the determinant representation of the partition function discovered by
Izergin \cite{izergin}.

This method is more flexible than Zeilberger's original proof and Kuperberg also used it in \cite{kup} to obtain many enumeration or equinumeration results for various symmetry classes of ASM's, most of them having been conjectured by Robbins \cite{robbins}. Among these results can be found the following remarkable one.

\begin{theo}\label{theo:kup}
{\em (Kuperberg).}
The number $A_{\textsc{QT}}(4N)$ of ASM's of size $4N$ invariant under a quarter-turn (QTASM's) is related to the number $A(N)$ of (unrestricted) ASM's of size $N$ and to the number $A_{\textsc{HT}}(2N)$ of ASM's of size $2N$ invariant under a half-turn (HTASM's) by the formula:
\begin{equation}\label{eq:kup}
A_{\textsc{QT}}(4N)=A_{\textsc{HT}}(2N) A(N)^2.
\end{equation}
\end{theo}

More recently, Razumov and Stroganov \cite{RS} applied Kuperberg's strategy to settle the following result relative to QTASM's of odd size, also conjectured by Robbins \cite{robbins} .

\begin{theo}\label{theo:RS}
{\em (Razumov, Stroganov).}
The numbers of QTASM's of odd size are given by the following formulas, where $A_{\textsc{HT}}(2N+1)$ is the number of HTASM's of size $2N+1$:
  \begin{align}
    A_{\textsc{QT}}(4N-1) & =  A_{\textsc{HT}}(2N-1) A(N)^2\label{eq:QT_m1}\\
    A_{\textsc{QT}}(4N+1) & =  A_{\textsc{HT}}(2N+1) A(N)^2\label{eq:QT_p1}.
  \end{align}
\end{theo}

It is easy to observe (and will be proved in Section \ref{sec:qQTASM}) that the set of QTASM's of size $4N+2$ is empty. But, by slightly relaxing the symmetry condition at the center of the matrix, Duchon introduced in \cite{duchon,duchonHDR} the notion of ASM's quasi-invariant under a quarter turn (the definition will be given in Section \ref{sec:qQTASM}) whose class is non-empty in size $4N+2$. Moreover, he conjectured for these qQTASM's an enumeration formula that perfectly completes the three previous enumeration results on QTASM. It is the aim of this paper to establish this formula.

\begin{theo}\label{theo:form}
The number $A_{\textsc{QT}}(4N+2)$ of qQTASM of size $4N+2$ is given by:
\begin{equation}\label{eq:form}
A_{\textsc{QT}}(4N+2)=A_{\textsc{HT}}(2N+1) A(N)A(N+1).
\end{equation}
\end{theo}

This paper is organized as follows: in Section \ref{sec:qQTASM}, we define qQTASM's; in Section \ref{sec:icemodel}, we recall the definitions of square ice models, precise the parameters and the partition functions that we shall study, and give the formula corresponding to equation \pref{eq:form} at the level of partition functions; Section \ref{sec:proofs} is devoted to the proofs; open questions are presented in Section \ref{sec:open}.

\section{ASM's quasi-invariant under a quarter-turn}\label{sec:qQTASM}

The class of ASM's invariant under a rotation by a quarter-turn (QTASM) is non-empty in size $4N-1$, $4N$, and $4N+1$. But this is not the case in size $4N+2$.

\begin{lemm} 
There is no QTASM of size $4N+2$.
\end{lemm}
\proof
Let us suppose that $M$ is a QTASM of even size $2L$. Now we use the fact that the size of an ASM is given by the sum of its entries, and the symmetry of $M$ to write:
\begin{equation}
2L=\sum_{1\le i,j\le2L}M_{i,j}=4\times\sum_{1\le i,j\le L}M_{i,j}
\end{equation}
which implies that the size of $M$ has to be a multiple of $4$.
\endproof

Duchon introduced in \cite{duchon,duchonHDR} a notion of ASM's quasi-invariant under a quarter-turn, by slightly relaxing the symmetry condition at the center of the matrix. The definition is more simple when considering the height matrix associated to the ASM, but can also be given directly.

\begin{defi}
An ASM $M$ of size $4N+2$ is said to be {\em quasi-invariant under a quarter-turn} (qQTASM) if its entries satisfy the quarter-turn symmetry
\begin{equation}
M_{4N+3-j,4N+3-i}=M_{i,j}
\end{equation}
except for the four central entries $(M_{2N+1,2N+1},M_{2N+1,2N+2},M_{2N+2,2N+1},M_{2N+2,2N+2})$ that have to be either $(0,-1,-1,0)$ or $(1,0,0,1)$.
\end{defi}

We give below two examples of qQTASM's of size $6$, with the two possible patterns at the center.
$$\left(
\begin{array}{cccccc}
0&0&0&1&0&0\cr
0&0&1&0&0&0\cr
1&0&0&-1&1&0\cr
0&1&-1&0&0&1\cr
0&0&0&1&0&0\cr
0&0&1&0&0&0\cr
\end{array}
\right)
\ \ \ \ \ \ \ \ \ \ \ \ \ \ \ \ 
\left(
\begin{array}{cccccc}
0&0&1&0&0&0\cr
0&1&-1&0&1&0\cr
0&0&1&0&-1&1\cr
1&-1&0&1&0&0\cr
0&1&0&-1&1&0\cr
0&0&0&1&0&0\cr
\end{array}
\right)$$

In the next section, we associate square ice models to ASM's with various types of symmetry.

\section{Square ice models and partition functions}\label{sec:icemodel}

\subsection{Notations}

Using Kuperberg's method we introduce square ice models associated to ASM's, HTASM's and (q)QTASM's. We recall here the main definitions and refer to \cite{kup} for details and many examples. 

Let $a\in\mathbb{C}$ be a global parameter.
For any complex number $x$ different from zero, we denote $\overline{x}=1/x$, and we define:
\begin{equation}
  \sigma(x) =  x - \overline{x}.
\end{equation}

Let $G$ denote some graph\footnote{Actually, our ``graphs'' are planar
  graphs together with a plane embedding.} where every vertex has
degree $1$, $2$ or $4$, with a fixed orientation attached to each edge incident to a vertex of degree $1$. An {\em ice state} of $G$ is an orientation of
the remaining edges such that every tetravalent vertex has exactly two incoming
and two outgoing edges, and each vertex of degree $2$ has either two
incoming or two outgoing edges.

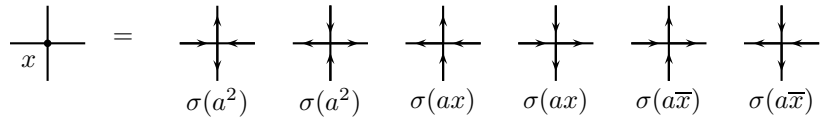
\begin{figure}[htbp]
  \begin{center}
\psset{unit=0.5cm}
\begin{pspicture}[.5](3,3)
\rput(1.5,1.5){
   \psline(-1,0)(1,0)
   \psline(0,-1)(0,1)
   \psdot[dotsize=.2]
   \rput(-.5,-.5){$x$}
  }
\end{pspicture}
$=$
    \begin{pspicture}[.5](19,5)
\rput(2,2.5){%
    \rput(-1,0){\flecheD}%
    \rput(1,0){\flecheG}%
    \rput(0,0){\flecheH}%
    \rput(0,0){\flecheB}%
    \rput[t](0,-1.2){$\sigma(a^2)$}%
}
\rput(5,2.5){%
    \rput(0,0){\flecheG}%
    \rput(0,0){\flecheD}%
    \rput(0,-1){\flecheH}%
    \rput(0,1){\flecheB}%
    \rput[t](0,-1.2){$\sigma(a^2)$}%
}
\rput(8,2.5){%
    \rput(0,0)\flecheG
    \rput(0,0)\flecheH
    \rput(1,0)\flecheG
    \rput(0,-1)\flecheH
    \rput[t](0,-1.2){$\sigma(ax)$}
}
\rput(11,2.5){
    \rput(0,0)\flecheD
    \rput(0,0)\flecheB
    \rput(-1,0)\flecheD
    \rput(0,1)\flecheB
    \rput[t](0,-1.2){$\sigma(ax)$}
}
\rput(14,2.5){
    \rput(-1,0)\flecheD
    \rput(0,0)\flecheD
    \rput(0,-1)\flecheH
    \rput(0,0)\flecheH
    \rput[t](0,-1.2){$\sigma(a\overline{x})$}
}
\rput(17,2.5){
    \rput(0,0)\flecheG
    \rput(1,0)\flecheG
    \rput(0,0)\flecheB
    \rput(0,1)\flecheB
    \rput[t](0,-1.2){$\sigma(a\overline{x})$}
}
    \end{pspicture}
  \end{center}
\caption{The 6 possible orientations and their associated weights}
\label{fig:poids_6V}
\end{figure}

A parameter $x\neq 0$ is assigned to one of the angles between
consecutive incident edges around each tetravalent vertex of the graph
$G$.  Then this vertex gets a weight, which depends on the orientation
of incident edges, as shown on Figure~\ref{fig:poids_6V}; the reader
can check that weights are unchanged if one of the angle parameters is
moved an adjacent angle of the same vertex, and simultaneously
replaced by its inverse.

It is sometimes easier to assign parameters, not to each vertex of the
graph, but to the lines that compose the graph. In this case, the
weight of a vertex is defined as:

\begin{displaymath}
\psset{unit=.5cm}
  \begin{pspicture}[.45](2,2)
\psline(1,0)(1,2)\psline(0,1)(2,1)
\rput[r](-.2,1){$x$}\rput[t](1,-.2){$y$}
  \end{pspicture}\ =\ 
  \begin{pspicture}[.45](2,2)
\psline(0,1)(2,1)\psline(1,0)(1,2)
\rput(.25,.25){$x\overline{y}$}
  \end{pspicture}
\end{displaymath}

\bigskip

When this convention is used, a parameter explicitly written at a
vertex replaces the quotient of the parameters of the lines.

We will put a dotted line to indicate that the parameter of a line is
different on the two sides of the dotted line.  

Vertices with degree $2$ do not get a parameter (they are only used to
force the two incident edges to have opposite orientations), and get
weight $1$.

\begin{displaymath}
  \psset{unit=.5cm}
  \begin{pspicture}[.6](2,1)
    \psline(0,.8)(2,.8)\rput(1,.8){\psdot[dotsize=.3]}
  \end{pspicture} = 
  \begin{pspicture}[.6](5,1)
    \rput(0,.8)\flecheD
    \rput(2,.8)\flecheG
    \rput[t](1,0.4){$1$}
    \rput(4,.8)\flecheG
    \rput(4,.8)\flecheD
    \rput[t](4,0.4){$1$}
    \rput(1,.8){\psdot[dotsize=.2]}
    \rput(4,.8){\psdot[dotsize=.2]}
  \end{pspicture}
\end{displaymath}


The partition function of a given ice graph is then defined as the sum,
over all its ice states, of the products of weights of all vertices.

To simplify notations, we will denote by $X_{N}$ the vector of variables
$(x_1,\dots, x_N)$. We use the notation
$X\backslash x$ to denote the vector $X$ without the variable $x$.

\subsection{Partition functions for classes of ASM's}

\begin{figure}[htbp]
  \begin{displaymath}
    Z(N;x_1,\dots, x_{N}, x_{N+1},\dots, x_{2N}) =
    \psset{unit=.5cm}
    \begin{pspicture}[.5](7,7)
      \rput(1,0){\IceSquare{6}}
      \rput[r](1,1){$x_{1}$}\rput[r](1,2){$x_{2}$}\rput[r](1,6){$x_{N}$}
      \rput[t](2,-.1){$x_{N+1}$}\rput[t](3,-.1){}\rput[t](7,-.1){$x_{2N}$}
    \end{pspicture}
  \end{displaymath}
\caption{Partition function for ASM's of size $N$}
\label{fig:Z}
\end{figure}
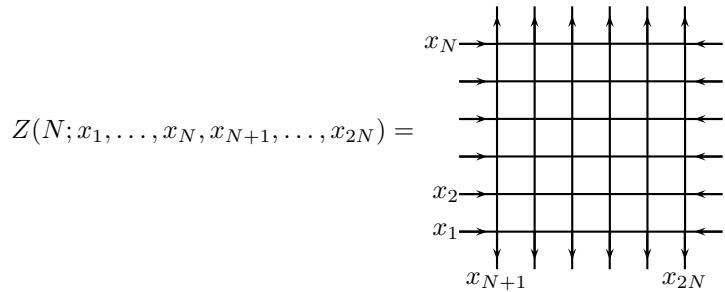

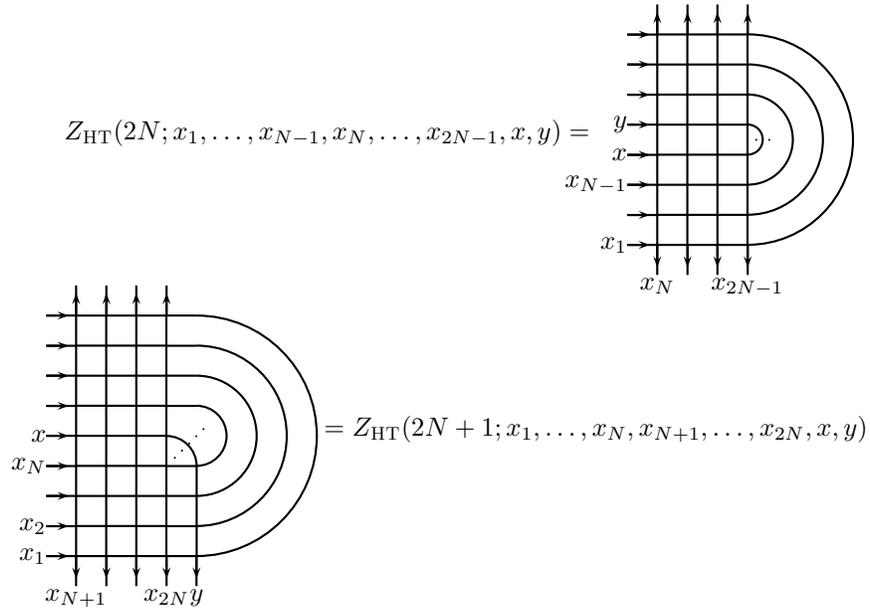
\begin{figure}[htbp]
  \begin{align*}
    Z_{\textsc{HT}}(2N; x_1,\dots,x_{N-1},x_{N},\dots,x_{2N-1},x,y)  = 
\psset{unit=.4cm}
    \begin{pspicture}[.5](9,9)
      \rput(1,0){\HTIceEven{4}}
      \rput[r](1,1){$x_{1}$}\rput[r](1,3){$x_{N-1}$}\rput[r](1,5){$y$}
      \rput[r](1,4){$x$}
      \rput[t](2,-.1){$x_{N}$}\rput[t](3,-.1){}\rput[t](5,-.1){$x_{2N-1}$}
    \end{pspicture}\\ 
\psset{unit=.4cm}
\begin{pspicture}[.5](10,10)
  \rput(1,0){\HTIceOdd{4}}
  \rput[r](1,1){$x_{1}$}\rput[r](1,2){$x_{2}$}\rput[r](1,4){$x_{N}$}
  \rput[r](1,5){$x$}
  \rput[t](2,-.1){$x_{N+1}$}\rput[t](3,-.1){}\rput[t](5,-.1){$x_{2N}$}
  \rput[t](6,-.1){$y$}
\end{pspicture}  =  
Z_{\textsc{HT}}(2N+1;x_1,\dots,x_{N},x_{N+1},\dots,x_{2N},x,y) 
  \end{align*}
\caption{Partition functions for HTASM's}
\label{fig:ZHT}
\end{figure}

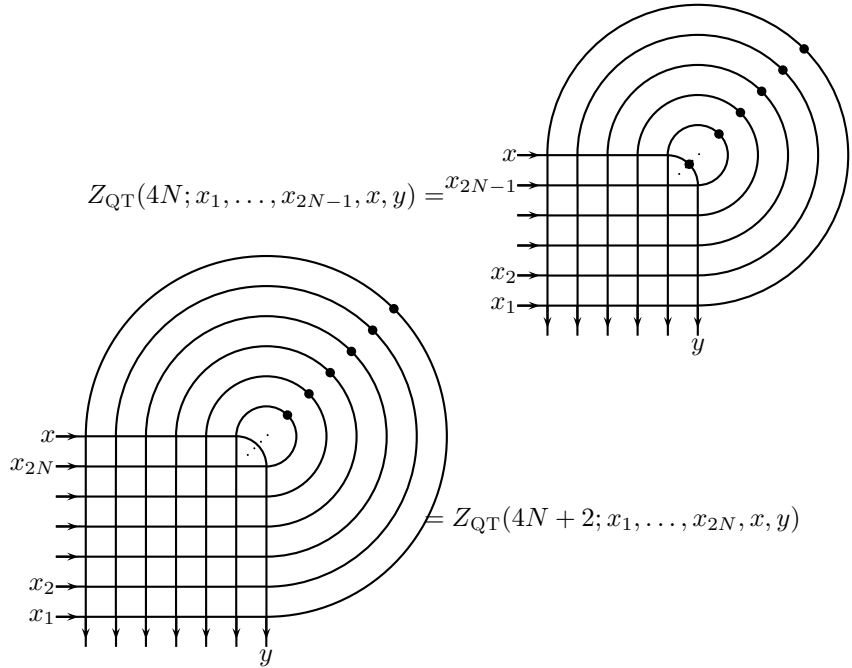
\begin{figure}[htbp]
\begin{align*}
  Z_{\textsc{QT}}(4N;x_1,\dots,x_{2N-1},x,y) & =\ \  \ \   
\psset{unit=.4cm}
  \begin{pspicture}[.4](11,11)
    \rput(1,0){\QTIce{6}}
    \rput[r](1,1){$x_1$}\rput[r](1,2){$x_2$}\rput[r](1,5){$x_{2N-1}$}
    \rput[r](1,6){$x$}\rput[t](7,-.1){$y$}
  \end{pspicture}\\
\psset{unit=.4cm}
  \begin{pspicture}[.4](12,10)
    \rput(0,0){\qQTIce{7}}
    \rput[r](0,1){$x_1$}\rput[r](0,2){$x_2$}\rput[r](0,6){$x_{2N}$}
    \rput[r](0,7){$x$}\rput[t](7,-.1){$y$}
  \end{pspicture} 
  & =  Z_{\textsc{QT}}(4N+2; x_1,\dots, x_{2N},x,y)
\end{align*}
\caption{Partition functions for (q)QTASM of even size}
\label{fig:ZQT}
\end{figure}

We give in Figures \ref{fig:Z}, \ref{fig:ZHT}, and \ref{fig:ZQT} the
ice models corresponding to the classes of ASM's that we shall study,
and their partition functions. The bijection between (unrestricted)
ASM's and states of the square ice model with ``domain wall boundary''
is now well-known (\cf \cite{kup}), and the bijections for the other
symmetry classes may be easily checked in the same way. The
correspondence between orientations of the ice model and entries of
ASM's is given in Figure \ref{fig:bij}.

\begin{figure}[htbp]
  \begin{center}
\psset{unit=0.5cm}
    \begin{pspicture}[.5](19,5)
\rput(2,2.5){%
    \rput(-1,0){\flecheD}%
    \rput(1,0){\flecheG}%
    \rput(0,0){\flecheH}%
    \rput(0,0){\flecheB}%
    \rput[t](0,-2.2){$1$}%
}
\rput(5,2.5){%
    \rput(0,0){\flecheG}%
    \rput(0,0){\flecheD}%
    \rput(0,-1){\flecheH}%
    \rput(0,1){\flecheB}%
    \rput[t](0,-2.2){$-1$}%
}
\rput(8,2.5){%
    \rput(0,0)\flecheG
    \rput(0,0)\flecheH
    \rput(1,0)\flecheG
    \rput(0,-1)\flecheH
    \rput[t](0,-2.2){$0$}%
}
\rput(11,2.5){
    \rput(0,0)\flecheD
    \rput(0,0)\flecheB
    \rput(-1,0)\flecheD
    \rput(0,1)\flecheB
    \rput[t](0,-2.2){$0$}%
}
\rput(14,2.5){
    \rput(-1,0)\flecheD
    \rput(0,0)\flecheD
    \rput(0,-1)\flecheH
    \rput(0,0)\flecheH
    \rput[t](0,-2.2){$0$}%
}
\rput(17,2.5){
    \rput(0,0)\flecheG
    \rput(1,0)\flecheG
    \rput(0,0)\flecheB
    \rput(0,1)\flecheB
    \rput[t](0,-2.2){$0$}%
}
    \end{pspicture}
  \end{center}
\caption{The correspondence between ice states and ASM's}
\label{fig:bij}
\end{figure}
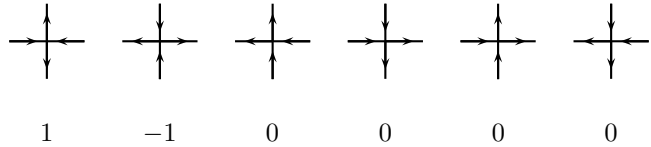

The reader may notice that the grid used to define
$Z_{\textsc{QT}}(4N)$ sligthly differs from the one used by Kuperberg
(the central vertices are treated in a different manner, and the line
$xy$ only carries a single parameter $x_{2N}$ in Kuperberg's model).
$Z_{\textsc{HT}}(2N)$ also appears in Kuperberg's paper \cite{kup} with a single
parameter on the $xy$ line; $Z_{\textsc{HT}}(2N+1)$ appears
identically in Razumov and Stroganov's paper~\cite{RSHT} (where a
different convention is used for the weights of vertices).

With these notations, Theorem \ref{theo:form} will be a consequence of the following one which addresses the concerned partition functions.

\begin{theo}\label{theo:main}
When $a=\omega_6=\exp(i\pi/3)$, one has for $N\ge 1$:
\begin{equation}\label{eq:main1}
\si(a)Z_{\textsc{QT}}(4N;X_{2N-1},x,y)=Z_{\textsc{HT}}(2N;X_{2N-1},x,y)Z(N;X_{2N-1},x)Z(N;X_{2N-1},y)
\end{equation}
and
\begin{equation}\label{eq:main2}
\si(a)Z_{\textsc{QT}}(4N+2;X_{2N},x,y)=Z_{\textsc{HT}}(2N+1;X_{2N},x,y)Z(N;X_{2N})Z(N+1;X_{2N},x,y).
\end{equation}
\end{theo}

Equation \pref{eq:main2} is new; Equation \pref{eq:main1} is due
to Kuperberg \cite{kup} for the case $x=y$.  To see that Theorem
\ref{theo:main} implies Theorem \ref{theo:form} (and Theorem
\ref{theo:kup}), we just have to observe that when $a=\omega_6$ and
all the variables are set to $1$, then the weight at each vertex is
$\si(a)=\si(a^2)=i\sqrt{3}$ thus the partition function reduces (up to
multiplication by $\si(a)^{\rm number\ of\ vertices}$) to the number
of states. This is summarized in the following proposition, where
$\mathbf{1}$ denotes the vector of all variables set to $1$.

\begin{prop}
  For $a=e^{i\pi/3}$, we have:
  \begin{align}
    Z(N;\mathbf{1}) & =  (i\sqrt{3})^{N^2} A(N)\\
    Z_{\textsc{HT}}(2N;\mathbf{1}) & =  (-1)^{N} 3^{N^2} A_{\textsc{HT}}(2N)\\
    Z_{\textsc{HT}}(2N+1;\mathbf{1}) & =  3^{N^2+N} A_{\textsc{HT}}(2N+1)\\
    Z_{\textsc{QT}}(4N;\mathbf{1}) & =  -i. 3^{2N^2-1/2} A_{\textsc{QT}}(4N)\\
    Z_{\textsc{QT}}(4N+2;\mathbf{1}) & =  3^{2N^2+2N} A_{\textsc{QT}}(4N+2).
  \end{align}
\end{prop}

\section{Proofs}\label{sec:proofs}

To prove Theorem \ref{theo:main}, the method, inspired from \cite{kup}, is to identify both sides of equations \pref{eq:main1} and \pref{eq:main2} as Laurent polynomials, and to produce as many specializations of the variables that verify the equalities, as needed to imply these equations in full generality.

In previous works \cite{kup,RS}, the final point in proofs is the evaluation of determinants or Pfaffians; in our proof of Theorem~\ref{theo:main}, 
we are able to avoid this computation by using symmetry properties.

\subsection{Laurent polynomials}

Since the weight of any vertex is a Laurent polynomial in the variables $x_i$, $x$ and $y$, the partition functions are Laurent polynomials in these variables. Moreover they are centered Laurent polynomials, \ie their lowest degree is the negative of their highest degree (called the half-width of the polynomial). In order to divide by two the number of non-zero coefficients (hence the number of required specializations) in $x$, we shall deal with Laurent polynomials of given parity in this variable. To do so, we group together the states with a given orientation (indicated as subscripts in the following notations) at the edge where the parameters $x$ and $y$ meet.

So let us consider the partition functions: 
\begin{itemize}
\item $\ZqConv(4N;X_{2N-1},x,y)$ and
$\ZqDiv(4N;X_{2N-1},x,y)$, respectively odd and even parts of
$Z_{\textsc{QT}}(4N;X_{2N-1},x,y)$ in $x$;
\item $\ZqUpleft(4N+2;X_{2N},x,y)$ and
$\ZqDownright(4N+2;X_{2N},x,y)$, respectively odd and even parts of
$Z_{\textsc{QT}}(4N+2;X_{2N},x,y)$in $x$;
\item $\ZhUp(2N;X_{2N-1},x,y)$ and
$\ZhDown(2N;X_{2N-1},x,y)$, respectively parts with the parity of $N$ and of $N-1$ of
$Z_{\textsc{HT}}(2N;X_{2N-1},x,y)$ in $x$; 
\item and
$\ZhUpleft(2N+1;X_{2N},x,y)$ and
$\ZhDownright(2N+1;X_{2N},x,y)$, respectively parts with the parity of $N-1$ and of $N$ of
$Z_{\textsc{HT}}(2N+1;X_{2N},x,y)$ in $x$.
\end{itemize}

With these notations, Equations \pref{eq:main1} and \pref{eq:main2}
are equivalent to the following:
\begin{align}
 \!\! \sigma(a)\ZqConv(4N;X_{2N-1},x,y) &\!\! =\!\! 
  \ZhUp(2N;X_{2N-1},x,y) Z(N;X_{2N-1},x)
  Z(N;X_{2N-1},y), \label{eq:main1.1}\\
 \!\! \sigma(a)\ZqDiv(4N;X_{2N-1},x,y) &\!\! =\!\! 
  \ZhDown(2N;X_{2N-1},x,y) Z(N;X_{2N-1},x)
  Z(N;X_{2N-1},y), \label{eq:main1.2}\\
 \!\!  \sigma(a)\ZqDownright(4N+2;X_{2N},x,y) &\!\! =\!\! 
  \ZhDownright(2N+1;X_{2N},x,y)
  Z(N+1;X_{2N},x,y) Z(N;X_{2N}), \label{eq:main2.1}\\
 \!\! \sigma(a)\ZqUpleft(4N+2;X_{2N},x,y) &\!\! =\!\! 
  \ZhUpleft(2N+1;X_{2N},x,y)
  Z(N+1;X_{2N},x,y) Z(N;X_{2N}). \label{eq:main2.2}
\end{align}

\begin{lemm}\label{lemm:Laurent}
  Both left-hand side and right-hand side of Equations
  (\ref{eq:main1.1}-\ref{eq:main2.2}) are centered Laurent polynomials
  in the variable $x$, odd or even, of respective half-widths $2N-1$,
  $2N-2$, $2N$, and $2N-1$. Thus, to prove each of these identities it
  is sufficient to exhibit specializations of $x$ for which the
  equality is true, and in number strictly exceeding the half-width.
\end{lemm}
\proof
To compute the half-width of these partition functions, we have to count the number of vertices in the ice models, and take note that non-zero entries of the ASM (\ie the first two orientations of Figure  \ref{fig:poids_6V}) give constant weight $\si(a^2)$. Also, a line whose orientation changes (respectively does not change) between endpoints must have an odd (respectively even) number of these $\pm 1$ entries.

We give the details for Equation \pref{eq:main1.1}:
\begin{itemize}
\item The term $Z(N;X_{2N-1},y)$ is a constant in $x$.
\item For $Z(N;X_{2N-1},x)$, the variable $x$ appears in the parameter 
of the $N$ vertices of the rightmost vertical line. On this line,
for each state of the model,
exactly one vertex gives a constant weight $\sigma(a^2)$, 
the other $N-1$ contribute for $1$ to the half-width. 
\item For $\ZhUp(2N;X_{2N-1},x,y)$, we have $N$ vertices on the line 
that carries the parameter $x$,
and an even number of them gives a constant weight.
\item For $\ZqConv(4N;X_{2N-1},x,y)$, we have in the same manner
$2N-1$ vertices that carries the parameter $x$, and an even number of them
gives a constant weight $\sigma(a^2)$.
\end{itemize}
This proves that both the left-hand side and the right-hand side of
equation \pref{eq:main1.1} are odd Laurent plynomial of half-width
$2N-1$. The assertions on Equations (\ref{eq:main1.2}-\ref{eq:main2.2})
are treated in the same way.
\endproof

\subsection{Symmetries}

To produce many specializations from one, we shall use symmetry properties of the partition functions. The crucial tool to prove this is the Yang-Baxter equation that we recall below.

\begin{lemm}\label{lemm:YB}{\em [Yang-Baxter equation]}
  If $xyz=\overline{a}$, then
  \begin{equation}
    \begin{pspicture}[.5](2,2)
      \SpecialCoor
      \rput(1,1){
	\psarc(-1.732,0){2}{330}{30}
        \psarc(.866,1.5){2}{210}{270}
        \psarc(.866,-1.5){2}{90}{150}
	\rput(1;60){$x$}
        \rput(1;180){$y$}
        \rput(1;300){$z$}
      }
    \end{pspicture}
     = 
     \begin{pspicture}[.5](2,2)
      \rput(1,1){
	\SpecialCoor
	\psarc(1.732,0){2}{150}{210}
        \psarc(-.866,1.5){2}{270}{330}
        \psarc(-.866,-1.5){2}{30}{90}
        \rput(1;240){$x$}
        \rput(1;0){$y$}
        \rput(1;120){$z$}
      }
     \end{pspicture}.
\label{eq:YB}
  \end{equation}
\end{lemm}

The following lemma gives a (now classical) example of use of the Yang-Baxter equation.

\begin{lemm}\label{lemm:echange_lignes}
  \begin{equation}
    \psset{unit=.5cm}
    \begin{pspicture}[.4](5.5,2)
      \rput[r](1,.5){$x$}\rput[r](1,1.5){$y$}
      \rput(1,.5){\flecheD}\rput(1,1.5){\flecheD}
      \psline(2,0)(2,2)\psline(3,0)(3,2)\psline(4.5,0)(4.5,2)
      \psline(2,.5)(3.5,.5)\psline(2,1.5)(3.5,1.5)
      \psline(4,.5)(4.5,.5)\psline(4,1.5)(4.5,1.5)
      \rput(3.75,1){$\dots$}
      \rput(5.5,.5){\flecheG}
      \rput(5.5,1.5){\flecheG}
    \end{pspicture} = 
    \begin{pspicture}[.4](5.5,2)
      \rput[r](1,.5){$y$}\rput[r](1,1.5){$x$}
      \rput(1,.5){\flecheD}\rput(1,1.5){\flecheD}
      \psline(2,0)(2,2)\psline(3,0)(3,2)\psline(4.5,0)(4.5,2)
      \psline(2,.5)(3.5,.5)\psline(2,1.5)(3.5,1.5)
      \psline(4,.5)(4.5,.5)\psline(4,1.5)(4.5,1.5)
      \rput(3.75,1){$\dots$}
      \rput(5.5,.5){\flecheG}
      \rput(5.5,1.5){\flecheG}      
    \end{pspicture}.
    \label{eq:symetrie_Z}
  \end{equation}
\end{lemm}

\proof
We multiply the left-hand side by $\sigma(a\overline{z})$, with
  $z=\overline{a}x\overline{y}$. We get
  \begin{align*} \psset{unit=.5cm}
    \sigma(a\overline{z})
    \begin{pspicture}[.4](5.5,2)
      \rput[r](1,.5){$x$}\rput[r](1,1.5){$y$}
      \rput(1,.5){\flecheD}\rput(1,1.5){\flecheD}
      \psline(2,0)(2,2)\psline(3,0)(3,2)\psline(4.5,0)(4.5,2)
      \psline(2,.5)(3.5,.5)\psline(2,1.5)(3.5,1.5)
      \psline(4,.5)(4.5,.5)\psline(4,1.5)(4.5,1.5)
      \rput(3.75,1){$\dots$}
      \rput(5.5,.5){\flecheG}
      \rput(5.5,1.5){\flecheG}
    \end{pspicture} & = \psset{unit=.5cm}
    \begin{pspicture}[.4](6.5,2)
      \rput[r](1,.5){$y$}\rput[r](1,1.5){$x$}
      \rput(1,.5){\flecheD}\rput(1,1.5){\flecheD}
      \rput(2,.5){\Hcrossing}\rput[r](2.4,1){$z$}
      \rput(1,0){
	\psline(2,0)(2,2)\psline(3,0)(3,2)\psline(4.5,0)(4.5,2)
	\psline(2,.5)(3.5,.5)\psline(2,1.5)(3.5,1.5)
	\psline(4,.5)(4.5,.5)\psline(4,1.5)(4.5,1.5)
	\rput(3.75,1){$\dots$}
	\rput(5.5,.5){\flecheG}
	\rput(5.5,1.5){\flecheG}}
    \end{pspicture} \\ 
      & =
    \psset{unit=.5cm}
    \begin{pspicture}[.4](5.5,2)
      \rput[r](1,.5){$y$}\rput[r](1,1.5){$x$}
      \rput(1,.5){\flecheD}\rput(1,1.5){\flecheD}
      \psline(2,0)(2,2)
      \rput(2,.5){\Hcrossing}\rput[r](2.4,1){$z$}
      \psline(3,0)(3,2)\psline(4.5,0)(4.5,2)
      \psline(3,.5)(3.5,.5)\psline(3,1.5)(3.5,1.5)
      \rput(3.75,1){$\dots$}
      \psline(4,.5)(4.5,.5)\psline(4,1.5)(4.5,1.5)
      \rput(5.5,.5){\flecheG}
      \rput(5.5,1.5){\flecheG}
    \end{pspicture} \\
      & = 
    \psset{unit=.5cm}
    \begin{pspicture}[.4](6.5,2)
      \rput[r](1,.5){$y$}\rput[r](1,1.5){$x$}
      \rput(1,.5){\flecheD}\rput(1,1.5){\flecheD}
      \psline(2,0)(2,2)
      \psline(3,0)(3,2)
      \psline(2,.5)(3.5,.5)\psline(2,1.5)(3.5,1.5)
      \rput(3.75,1){$\dots$}
      \psline(4,.5)(4.5,.5)\psline(4,1.5)(4.5,1.5)
      \rput(4.5,.5){\Hcrossing}\rput[r](4.9,1){$z$}
      \psline(5.5,0)(5.5,2)
      \rput(6.5,.5){\flecheG}
      \rput(6.5,1.5){\flecheG}
    \end{pspicture} \\
      & = 
    \psset{unit=.5cm}
    \begin{pspicture}[.4](6.5,2)
      \rput[r](1,.5){$y$}\rput[r](1,1.5){$x$}
      \rput(1,.5){\flecheD}\rput(1,1.5){\flecheD}
      \psline(2,0)(2,2)\psline(3,0)(3,2)\psline(4.5,0)(4.5,2)
      \psline(2,.5)(3.5,.5)\psline(2,1.5)(3.5,1.5)
      \psline(4,.5)(4.5,.5)\psline(4,1.5)(4.5,1.5)
      \rput(3.75,1){$\dots$}
      \rput(4.5,.5){\Hcrossing}\rput[l](5.1,1){$z$}
      \rput(6.5,.5){\flecheG}
      \rput(6.5,1.5){\flecheG}      
    \end{pspicture} \\
      & = 
    \psset{unit=.5cm}
    \begin{pspicture}[.4](5.5,2)
      \rput[r](1,.5){$y$}\rput[r](1,1.5){$x$}
      \rput(1,.5){\flecheD}\rput(1,1.5){\flecheD}
      \psline(2,0)(2,2)\psline(3,0)(3,2)\psline(4.5,0)(4.5,2)
      \psline(2,.5)(3.5,.5)\psline(2,1.5)(3.5,1.5)
      \psline(4,.5)(4.5,.5)\psline(4,1.5)(4.5,1.5)
      \rput(3.75,1){$\dots$}
      \rput(5.5,.5){\flecheG}
      \rput(5.5,1.5){\flecheG}      
    \end{pspicture} \sigma(a\overline{z})
  \end{align*}
\endproof

The same method, together with the easy transformation
  \begin{equation}
    \psset{unit=.5cm}
    \begin{pspicture}[.45](2.5,1)
      \Hcrossing
      \psarc(1,.5){.5}{270}{90}
      \rput[r](.3,.5){$z$}
    \end{pspicture} = \left(\sigma(az)+\sigma(a^2)\right)
    \left(
    \begin{pspicture}[.45](1,1)
      \rput(0,1){\flecheD}\rput(1,0){\flecheG}
    \end{pspicture}\ +\ 
    \begin{pspicture}[.45](1,1)
      \rput(0,0){\flecheD}\rput(1,1){\flecheG}
    \end{pspicture}\right)
    \label{eq:boucle}
  \end{equation}
gives the following lemma.

\begin{lemm}\label{lemm:echange_boucle}
  \begin{align}
    \psset{unit=.5cm}
    \begin{pspicture}[.4](5.7,2)
      \rput[r](1,.5){$x$}\rput[r](1,1.5){$y$}
      \rput(1,.5){\flecheD}\rput(1,1.5){\flecheD}
      \psline(2,0)(2,2)\psline(3,0)(3,2)\psline(4.5,0)(4.5,2)
      \psline(2,.5)(3.5,.5)\psline(2,1.5)(3.5,1.5)
      \psline(4,.5)(4.5,.5)\psline(4,1.5)(4.5,1.5)
      \psarc(4.5,1){.5}{270}{90}
      \psline[linestyle=dotted](4.7,1)(5.7,1)
      \rput(3.75,1){$\dots$}
    \end{pspicture} &=
    \psset{unit=.5cm}
    \frac{\sigma(a^2)+\sigma(x\overline{y})}{\sigma(a^2y\overline{x})}
    \begin{pspicture}[.4](5.7,2)
      \rput[r](1,.5){$y$}\rput[r](1,1.5){$x$}
      \rput(1,.5){\flecheD}\rput(1,1.5){\flecheD}
      \psline(2,0)(2,2)\psline(3,0)(3,2)\psline(4.5,0)(4.5,2)
      \psline(2,.5)(3.5,.5)\psline(2,1.5)(3.5,1.5)
      \psline(4,.5)(4.5,.5)\psline(4,1.5)(4.5,1.5)
      \psarc(4.5,1){.5}{270}{90}
      \psline[linestyle=dotted](4.7,1)(5.7,1)
      \rput(3.75,1){$\dots$}
    \end{pspicture} \label{eq:echange_boucle}\\
    \psset{unit=.5cm}
    \begin{pspicture}[.4](5.7,2)
      \rput[r](1,.5){$x$}\rput[r](1,1.5){$y$}
      \rput(1,.5){\flecheD}\rput(1,1.5){\flecheD}
      \psline(2,0)(2,2)\psline(3,0)(3,2)\psline(4.5,0)(4.5,2)
      \psline(2,.5)(3.5,.5)\psline(2,1.5)(3.5,1.5)
      \psline(4,.5)(4.5,.5)\psline(4,1.5)(4.5,1.5)
      \rput(4.5,1.5){\flecheD}
      \rput(5.5,.5){\flecheG}
      \rput(3.75,1){$\dots$}
    \end{pspicture} &=
    \psset{unit=.5cm}
\frac{\sigma(x\overline{y})}{\sigma(a^2y\overline{x})}
    \begin{pspicture}[.4](5.7,2)
      \rput[r](1,.5){$y$}\rput[r](1,1.5){$x$}
      \rput(1,.5){\flecheD}\rput(1,1.5){\flecheD}
      \psline(2,0)(2,2)\psline(3,0)(3,2)\psline(4.5,0)(4.5,2)
      \psline(2,.5)(3.5,.5)\psline(2,1.5)(3.5,1.5)
      \psline(4,.5)(4.5,.5)\psline(4,1.5)(4.5,1.5)
      \rput(5.5,1.5){\flecheG}
      \rput(4.5,.5){\flecheD}
      \rput(3.75,1){$\dots$}
    \end{pspicture}  + 
\frac{\sigma(a^2)}{\sigma(a^2y\overline{x})}
    \begin{pspicture}[.4](5.7,2)
      \rput[r](1,.5){$y$}\rput[r](1,1.5){$x$}
      \rput(1,.5){\flecheD}\rput(1,1.5){\flecheD}
      \psline(2,0)(2,2)\psline(3,0)(3,2)\psline(4.5,0)(4.5,2)
      \psline(2,.5)(3.5,.5)\psline(2,1.5)(3.5,1.5)
      \psline(4,.5)(4.5,.5)\psline(4,1.5)(4.5,1.5)
      \rput(4.5,1.5){\flecheD}
      \rput(5.5,.5){\flecheG}
      \rput(3.75,1){$\dots$}
    \end{pspicture} \label{eq:echange_boucle_a}\\
    \psset{unit=.5cm}
    \begin{pspicture}[.4](5.7,2)
      \rput[r](1,.5){$x$}\rput[r](1,1.5){$y$}
      \rput(1,.5){\flecheD}\rput(1,1.5){\flecheD}
      \psline(2,0)(2,2)\psline(3,0)(3,2)\psline(4.5,0)(4.5,2)
      \psline(2,.5)(3.5,.5)\psline(2,1.5)(3.5,1.5)
      \psline(4,.5)(4.5,.5)\psline(4,1.5)(4.5,1.5)
      \rput(4.5,.5){\flecheD}
      \rput(5.5,1.5){\flecheG}
      \rput(3.75,1){$\dots$}
    \end{pspicture} &=
    \psset{unit=.5cm}
\frac{\sigma(x\overline{y})}{\sigma(a^2y\overline{x})}
    \begin{pspicture}[.4](5.7,2)
      \rput[r](1,.5){$y$}\rput[r](1,1.5){$x$}
      \rput(1,.5){\flecheD}\rput(1,1.5){\flecheD}
      \psline(2,0)(2,2)\psline(3,0)(3,2)\psline(4.5,0)(4.5,2)
      \psline(2,.5)(3.5,.5)\psline(2,1.5)(3.5,1.5)
      \psline(4,.5)(4.5,.5)\psline(4,1.5)(4.5,1.5)
      \rput(5.5,.5){\flecheG}
      \rput(4.5,1.5){\flecheD}
      \rput(3.75,1){$\dots$}
    \end{pspicture}  + 
\frac{\sigma(a^2)}{\sigma(a^2y\overline{x})}
    \begin{pspicture}[.4](5.7,2)
      \rput[r](1,.5){$y$}\rput[r](1,1.5){$x$}
      \rput(1,.5){\flecheD}\rput(1,1.5){\flecheD}
      \psline(2,0)(2,2)\psline(3,0)(3,2)\psline(4.5,0)(4.5,2)
      \psline(2,.5)(3.5,.5)\psline(2,1.5)(3.5,1.5)
      \psline(4,.5)(4.5,.5)\psline(4,1.5)(4.5,1.5)
      \rput(4.5,.5){\flecheD}
      \rput(5.5,1.5){\flecheG}
      \rput(3.75,1){$\dots$}
    \end{pspicture} \label{eq:echange_boucle_b}
  \end{align}
\end{lemm}

We use Lemmas \ref{lemm:echange_lignes} and \ref{lemm:echange_boucle} to obtain symmetry properties of the partition functions, that we summarize below, where $m$ denotes either $2N$ or $2N+1$. 

\begin{lemm}\label{lemm:sym}
  The functions $Z(N;X_{2N})$ and 
  $Z_{\textsc{HT}}(2N+1;X_{2N},x,y)$ are symmetric separately in the two sets of variables $\{x_i,\ i\le N\}$ and $\{x_i,\ i\ge N+1\}$, the function 
  $Z_{\textsc{HT}}(2N;X_{2N-1},x,y)$ is symmetric separately in the two sets of variables $\{x_i,\ i\le N-1\}$ and $\{x_i,\ i\ge N\}$, 
  and the functions 
  $Z_{\textsc{QT}}(2m;X_{N-1},x,y)$ are symmetric in their variables $x_i$.

Moreover, $Z_{\textsc{QT}}(4N+2;\dots)$ is symmetric in its variables $x$ and
$y$, and we have a pseudo-symmetry for $Z_{\textsc{QT}}(4N;\dots)$ and $Z_{\textsc{HT}}(2N;\dots)$:
  \begin{align}
    Z_{\textsc{QT}}(4N;X_{2N-1},x,y) & =  
    \frac{\sigma(a^2)+\sigma(x\overline{y})}{\sigma(a^2y\overline{x})}
    Z_{\textsc{QT}}(4N;X_{2N-1},y,x),
    \label{eq:sym_QT4N_xy}\\
    Z_{\textsc{HT}}(2N;X_{2N-1},x,y) & =  
    \frac{\sigma(a^2)+\sigma(x\overline{y})}{\sigma(a^2y\overline{x})}
    Z_{\textsc{HT}}(2N;X_{2N-1},y,x).
    \label{eq:sym_HT2N_xy}
  \end{align}
\end{lemm}

\proof
  
For $Z(N;\dots)$, $Z_{\textsc{HT}}(m;\dots)$ and $Z_{\textsc{QT}}(2m;\dots)$, 
  the symmetry in two ``consecutive'' 
  variables $x_i$ and $x_{i+1}$ is a direct consequence of
  Lemma~\ref{lemm:echange_lignes}. 

For the (pseudo-)symmetry of $Z_{\textsc{QT}}(2m;\dots)$, we use the easy observations:
\begin{equation}\label{eq:2dots}
  \psset{unit=5mm}
    \begin{pspicture}[.4](2,2)
      \psline(1,0)(1,2)\psline(0,1)(2,1)
    \end{pspicture}
    \ =\  
     \begin{pspicture}[.4](2,2)
       \psline(1,0)(1,2)\psline(0,1)(2,1)
       \psdot[dotsize=.2](.5,1)
       \psdot[dotsize=.2](1,1.5)
       \psdot[dotsize=.2](1,.5)
       \psdot[dotsize=.2](1.5,1)
     \end{pspicture} 
\ \ \ \ {\rm and}\ \ \ \ 
     \begin{pspicture}[.4](2,2)
       \psline(0,1)(2,1)
     \end{pspicture}
    \ =\  
     \begin{pspicture}[.4](2,2)
       \psline(0,1)(2,1)
       \psdot[dotsize=.2](.5,1)
       \psdot[dotsize=.2](1.5,1)
     \end{pspicture}
\end{equation}
which  gives us the following modification of the grid in size $4N+2$:
  \begin{align*}
    Z_{\textsc{QT}}(4N+2;{X}_{2N},x,y) & = 
    \psset{unit=.3cm}
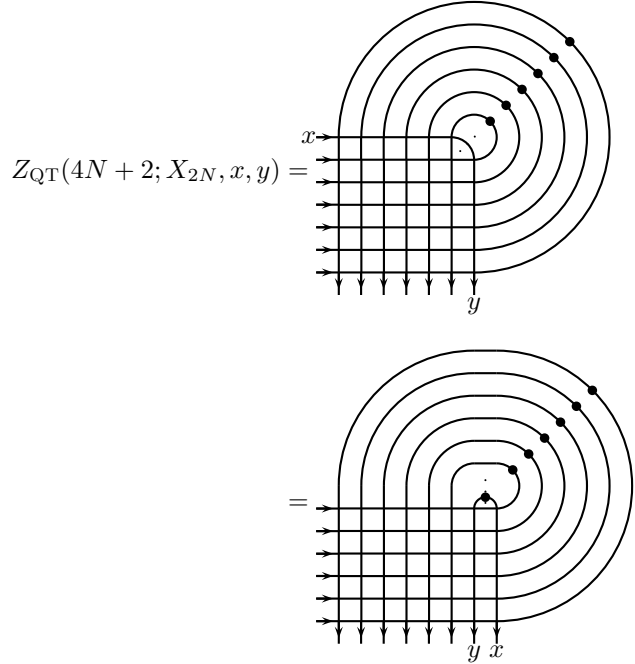
\begin{pspicture}[.4](12,13)
    \rput(0,0){\qQTIce{7}}
    \rput[r](0,7){$x$}\rput[t](7,-.1){$y$}
  \end{pspicture}\\
  & =  
\psset{unit=.3cm}
  \begin{pspicture}[.4](12,15)
    \rput(1,1){\IceGrid{8}{6}}
    \rput(0,1){\colD{6}}
    \rput(1,0){\linB{8}}
    \rput[t](7,-.1){$y$}
    \rput[t](8,-.1){$x$}
    \psarc(7.5,6){.5}{0}{180}
    \psdots[dotstyle=*](7.5,6.5)
    \psline[linestyle=dotted](7.5,6.2)(7.5,7.3)
    \multido{\i=1+1}{6}{
      \psline(\i,6)(\i,7)
      \psarc(7,7){\i}{90}{180}
      \psarc(8,7){\i}{270}{90}
      }
      \multido{\i=8+1}{6}{\psline(7,\i)(8,\i)}
    \rput(8,7){\SpecialCoor\multido{\i=1+1}{6}{\psdots[dotstyle=*](\i;45)}}
    \end{pspicture}\\
  \end{align*}

and in size $4N$:
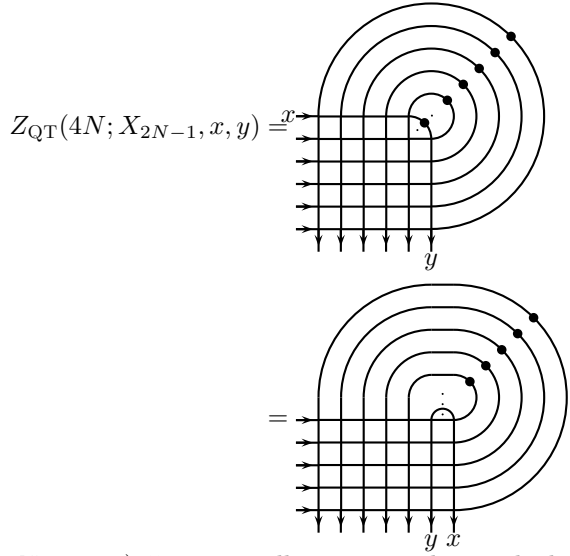
  \begin{align*}
Z_{\textsc{QT}}(4N;{X}_{2N-1},x,y) & =  
\psset{unit=.3cm}\begin{pspicture}[.4](10,13)
\rput(0,0){\QTIce{6}}
\rput[r](0,6){$x$}\rput[t](6,-.1){$y$}
\end{pspicture}\\
& = 
\psset{unit=.3cm}
  \begin{pspicture}[.4](12,12)
    \rput(1,1){\IceGrid{7}{5}}
    \rput(0,1){\colD{5}}
    \rput(1,0){\linB{7}}
    \rput[t](6,-.1){$y$}
    \rput[t](7,-.1){$x$}
    \psarc(6.5,5){.5}{0}{180}
    \psline[linestyle=dotted](6.5,5.2)(6.5,6.2)
    \multido{\i=1+1}{5}{
      \psline(\i,5)(\i,6)
      \psarc(6,6){\i}{90}{180}
      \psarc(7,6){\i}{270}{90}
      }
      \multido{\i=7+1}{5}{\psline(6,\i)(7,\i)}
    \rput(7,6){\SpecialCoor\multido{\i=1+1}{5}{\psdots[dotstyle=*](\i;45)}}
    \end{pspicture}
  \end{align*}

Now for $Z_{\textsc{QT}}(4N+2;\dots)$, Lemma~\ref{lemm:echange_lignes} allow us to exchange the lines
of parameters 
  $x$ and $y$, which proves the symmetry. for $Z_{\textsc{QT}}(4N;\dots)$,
  we apply Lemma~\ref{lemm:echange_boucle} to conclude.

The last assertion concerns $Z_{\textsc{HT}}(2N;\dots)$, and
  Lemma~\ref{lemm:echange_boucle} gives directly
  (\ref{eq:sym_HT2N_xy}) without any modification of the graph.

\endproof

\begin{rema}\label{rema:symoddeven}
It should be clear, but is useful to note, that we  have analogous properties for the even and odd parts of the partition functions.
\end{rema}

The next (and last) symmetry property, proved by Stroganov
\cite{IKdet}, appears when the parameter $a$ takes the special value
$\omega_6=\exp(i\pi/3)$. 

\begin{lemm}\label{lemm:symZomega6}
When $a=\omega_6=\exp(i\pi/3)$, the partition function $Z(N;X_{2N})$ is symmetric in {\bf all} its variables.
\end{lemm}

Stroganov proved this surprising symmetry property by a study of Izergin-Korepin determinant.
A proof only involving Yang-Baxter equation has recently been given in \cite{aval}.

\subsection{Specializations, recurrences}

The aim of this section is to give the value of the partition
functions in some specializations of the variable $x$ or $y$. The
first result is due to Kuperberg; the others are very similar.

\begin{lemm}{\em [specialization of $Z$; Kuperberg]}
\label{lemm:specZ}
If we denote
\begin{align*}
  A(x_{N+1},X_{2N}\backslash \{x_1,x_{N+1}\}) & = 
    \prod_{2\leq k\leq N} \sigma(a x_k \overline{x}_{N+1})
    \prod_{N+1\leq k\leq 2N} \sigma(a^2 x_{N+1} \overline{x}_k),\\
  \overline{A}(x_{N+1},X_{2N}\backslash \{x_1,x_{N+1}\}) & = 
    \prod_{2\leq k\leq N} \sigma(a x_{N+1} \overline{x}_k) 
    \prod_{N+1\leq k\leq 2N} \sigma(a^2 x_k \overline{x}_{N+1}),
\end{align*}
then we have:
  \begin{align}
 Z(N;{\bf \overline{a}x_{N+1}},X_{2N}\backslash x_1) & = 
\overline{A}(x_{N+1},X_{2N}\backslash \{x_1,x_{N+1}\}) Z(N-1;X_{2N}\backslash \{x_1,x_{N+1}\}), \label{eq:Zbax}\\
 Z(N;{\bf ax_{N+1}},X_{2N}\backslash x_1) & = 
A(x_{N+1},X_{2N}\backslash \{x_1,x_{N+1}\}) 
 Z(N-1; X_{2N}\backslash \{x_1,x_{N+1}\}). \label{eq:Zax}
  \end{align}
\end{lemm}

\proof
We recall the method to prove equation \pref{eq:Zbax}. 
We consider the crossing of the lines of parameter $x_1$ and $x_{N+1}$ which could of the following two types:

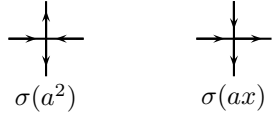
\begin{displaymath}
\psset{unit=.5cm}
    \begin{pspicture}[.5](19,5)
\rput(7,2.5){%
    \rput(-1,0){\flecheD}%
    \rput(1,0){\flecheG}%
    \rput(0,0){\flecheH}%
    \rput(0,0){\flecheB}%
    \rput[t](0,-1.2){$\sigma(a^2)$}%
}
\rput(12,2.5){
    \rput(0,0)\flecheD
    \rput(0,0)\flecheB
    \rput(-1,0)\flecheD
    \rput(0,1)\flecheB
    \rput[t](0,-1.2){$\sigma(ax)$}
}

    \end{pspicture}
\end{displaymath}

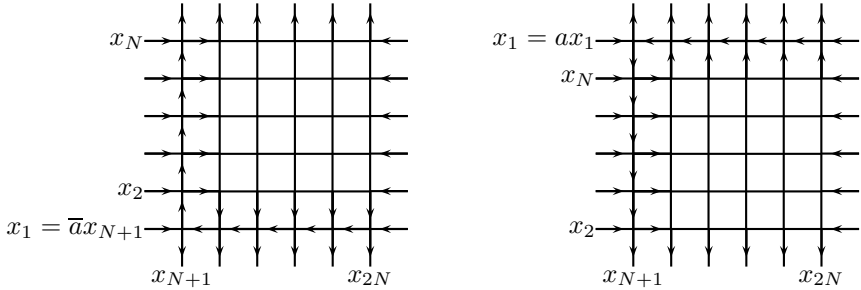
\begin{figure}[htbp]
    \begin{center}\psset{unit=.5cm}
      \begin{pspicture}(-1,-1)(19,8)
	\rput(0,0){\IceSquare{6}}
	\rput[r](0,1){$x_1=\overline{a}x_{N+1}$}\rput[r](0,6){$x_{N}$}
        \rput[r](0,2){$x_2$}
        \rput[t](1,-.1){$x_{N+1}$}\rput[t](6,-.1){$x_{2N}$}
	\multido{\i=1+1}{5}{\rput(1,\i){\flecheH}}
        \multido{\i=2+1}{5}{\rput(1,\i){\flecheD}\rput(\i,1){\flecheG}\rput(\i,2){\flecheB}}

	\rput(12,0){%
	  \IceSquare{6}
          \rput[r](0,6){$x_1=ax_{1}$}\rput[r](0,5){$x_{N}$}\rput[r](0,1){$x_2$}
          \rput[t](1,-.1){$x_{N+1}$}\rput[t](6,-.1){$x_{2N}$}
          \multido{\i=2+1}{5}{\rput(1,\i){\flecheB}\rput(\i,5){\flecheH}
             \rput(\i,6){\flecheG}}
          \multido{\i=1+1}{5}{\rput(1,\i){\flecheD}}
	}
      \end{pspicture}
      \caption{Fixed edges for (\ref{eq:Zbax}) on the left and
      (\ref{eq:Zax}) on the right}
    \label{fig:recZ}
    \end{center}
  \end{figure}

We observe that when $x_1=\bar a x_{N+1}$, the parameter of the vertex at the crossing of the two lines of parameter $x_1$ and $x_{N+1}$ is $x=x_1\bar x_{N+1}=\bar a$. 
Thus the weight of this vertex is $\si(a\ba)=\si(1)=0$ in the second orientation on the figure above. This forces the orientation of this vertex to be the first described in this figure.
But this orientation implies the orientation of all vertices in the row $x_1$ and in the column $x_{N+1}$, as shown on Figure \ref{fig:recZ}. 
The non-fixed part gives the partition function $Z$ in size $N-1$, without parameters $x_1$ and
  $x_{N+1}$, and the weights of the fixed part gives the factor $\overline{A}(\dots)$.

The case of \pref{eq:Zax} is similar, after using Lemma \ref{lemm:sym} to put the line $x_{N+1}$ at the top of the grid, as shown on Figure \ref{fig:recZ}.

\endproof

We will need the following application of the Yang-Baxter equation, which allows, under certain condition, a line with a change of parameter to go through a grid.

\begin{lemm}\label{lemm:tour_grille}
  \begin{equation}
    \psset{unit=.3cm}
    \begin{pspicture}[.5](9,8)
      \multido{\i=2+1}{6}{\psline(\i,0)(\i,8)}
      \multido{\i=2+1}{5}{\psline(0,\i)(9,\i)}
      \psline(1,7)(7,7)\psline(8,1)(8,6)
      \psarc(7,6){1}{0}{90}
      \psline[linestyle=dotted](7.2,6.2)(8.5,7.5)
      \rput[r](1,7){$x$}\rput[t](8,.9){$\overline{a}x$}
    \end{pspicture}
    = 
    \begin{pspicture}[.5](9,8)
      \multido{\i=2+1}{6}{\psline(\i,0)(\i,8)}
      \multido{\i=2+1}{5}{\psline(0,\i)(9,\i)}
      \psline(1,7)(1,2)\psline(2,1)(8,1)
      \psarc(2,2){1}{180}{270}
      \psline[linestyle=dotted](1.8,1.8)(.5,.5)
      \rput[l](8,1){$x$}\rput[b](1,7){$\overline{a}x$}
    \end{pspicture}
  \label{eq:tour_grille}
  \end{equation}
\end{lemm}

\proof
 We iteratively apply Lemma \ref{lemm:YB} on the rows, and row by row:
  \begin{align*}
    \psset{unit=.3cm}
    \begin{pspicture}[.4](10,4)
      \psline(0,2)(10,2)
      \multido{\i=2+2}{4}{\psline(\i,0)(\i,4)}
      \psline(1,3)(8,3)\psline(9,2)(9,1)\psarc(8,2){1}{0}{90}
      \psline[linestyle=dotted](8.2,2.2)(9.5,3.5)
      \rput[r](1,3){$x$}\rput[t](9,.9){$\overline{a}x$}
    \end{pspicture}
    & = 
    \psset{unit=.3cm}
    \begin{pspicture}[.4](10,4)
      \psline(0,2)(10,2)
      \multido{\i=2+2}{4}{\psline(\i,0)(\i,4)}
      \psline(1,3)(6,3)\psline(8,1)(9,1)\psarc(6,2){1}{0}{90}
      \psarc(8,2){1}{180}{270}
      \psline[linestyle=dotted](6.2,2.2)(7.5,3.5)
      \psline[linestyle=dotted](7.8,1.8)(6.5,.5)
      \rput[r](1,3){$x$}
      \rput[l](9,1){$x$}
      \rput[r](6.8,1.5){$\overline{a}x$}
    \end{pspicture}\\
    & = 
    \psset{unit=.3cm}
    \begin{pspicture}[.4](10,4)
      \psline(0,2)(10,2)
      \multido{\i=2+2}{4}{\psline(\i,0)(\i,4)}
      \psline(1,3)(4,3)\psline(6,1)(9,1)\psarc(4,2){1}{0}{90}
      \psarc(6,2){1}{180}{270}
      \psline[linestyle=dotted](4.2,2.2)(5.5,3.5)
      \psline[linestyle=dotted](5.8,1.8)(4.5,.5)
      \rput[r](1,3){$x$}
      \rput[l](9,1){$x$}
      \rput[r](4.8,1.5){$\overline{a}x$}
    \end{pspicture}\\
    & = 
    \psset{unit=.3cm}
    \begin{pspicture}[.4](10,4)
      \psline(0,2)(10,2)
      \multido{\i=2+2}{4}{\psline(\i,0)(\i,4)}
      \psline(1,3)(4,3)\psline(6,1)(9,1)\psarc(4,2){1}{0}{90}
      \psarc(6,2){1}{180}{270}
      \psline[linestyle=dotted](4.2,2.2)(5.5,3.5)
      \psline[linestyle=dotted](5.8,1.8)(4.5,.5)
      \rput[r](1,3){$x$}
      \rput[l](9,1){$x$}
      \rput[r](4.8,1.5){$\overline{a}x$}
    \end{pspicture}\\
    & = 
    \psset{unit=.3cm}
    \begin{pspicture}[.4](10,4)
      \psline(0,2)(10,2)
      \multido{\i=2+2}{4}{\psline(\i,0)(\i,4)}
      \psline(1,3)(1,2)\psline(2,1)(9,1)\psarc(2,2){1}{180}{270}
      \psline[linestyle=dotted](1.8,1.8)(.5,.5)
      \rput[b](1,3){$\overline{a}x$}\rput[l](9,1){$x$}
    \end{pspicture}.
  \end{align*}

\endproof

\begin{lemm}{\em [specialization of $Z_{\textsc{HT}}$]}
\label{lemm:specZHT}
If we denote
 \begin{align*}
   A_{H}^{1}(x_1,X_{2N}\backslash x_1) & = 
     \prod_{1\leq k\leq N} \sigma(a^2 x_1 \overline{x}_k)
     \prod_{N+1\leq k\leq 2N} \sigma(a x_k\overline{x}_1),\\
   \overline{A}_{H}^{1}(x_1,X_{2N}\backslash x_1) & = 
     \prod_{1\leq k\leq N} \sigma(a^2x_k\overline{x}_1)
     \prod_{N+1\leq k\leq 2N} \sigma(ax_1\overline{x}_k),\\
   A_{H}^{0}(x_{N},X_{2N-1}\backslash x_{N}) & =  
     \prod_{1\leq k\leq N-1}\sigma(a x_k \overline{x}_N)
     \prod_{N\leq k\leq 2N-1} \sigma(a^2 x_N \overline{x}_k),\\
   \overline{A}_{H}^{0}(x_{N},X_{2N-1}\backslash x_{N}) & = 
     \prod_{1\leq k\leq N-1} \sigma(ax_N\overline{x}_k)
     \prod_{N\leq k\leq 2N-1} \sigma(a^2x_k\overline{x}_N),
 \end{align*}
  then for $\star=\downright,\upleft,\upc,\downc$ and
  $\square=\downc,\upc,\downright,\upleft$ respectively, we have
  \begin{align}
\!\!\!\!\!\!\!\!Z_{\textsc{HT}}^{\star}(2N+1;X_{2N},x,\mathbf{ax_1}) &\!\!=\!\!
    A_{H}^{1}(x_1,X_{2N}\backslash x_{1}) 
    Z_{\textsc{HT}}^{\square}(2N;X_{2N}\backslash x_1, x_1,x),
    \label{equa:Zht_impair_ax}\\
\!\!\!\!\!\!\!\!Z_{\textsc{HT}}^{\square}(2N+1;X_{2N},x,{\bf \overline{a}x_1}) &\!\!=\!\!
    \overline{A}_{H}^{1}(x_1,X_{2N}\backslash x_1)
     Z_{\textsc{HT}}^{\star}(2N;X_{2N}\backslash x_1,x,x_1),
    \label{equa:Zht_impair_bax}\\
\!\!\!\!\!\!\!\!Z_{\textsc{HT}}^{\star}(2N;X_{2N-1},x,{\bf ax_N}) &\!\!=\!\!
  \sigma(ax\overline{x}_{N}) A_{H}^{0}(x_{N},X_{2N-1}\backslash x_{N})
  Z_{\textsc{HT}}^{\square}(2N-1;X_{2N-1}\backslash x_N,x,x_N),   
    \label{equa:Zht_pair_ax}\\
\!\!\!\!\!\!\!\!Z_{\textsc{HT}}^{\square}(2N;X_{2N-1},{\bf \overline{a}x_N},y) &\!\!=\!\!
  \sigma(ax_N\overline{y}) 
  \overline{A}_{H}^{0}(x_{N},X_{2N-1}\backslash x_{N})
   Z_{\textsc{HT}}^{\star}(2N-1;X_{2N-1}\backslash x_N,y,x_N).
    \label{equa:Zht_pair_bax}
  \end{align}
\end{lemm}

\proof

The proof is similar to the previous one, with the difference 
that before looking at fixed edges, we need to multiply the partition function by a given factor; 
we interpret this operation by a modification of the graph of the ice model, and apply 
Lemma~\ref{lemm:tour_grille}. It turns out that in each case, the additional factors 
are exactly cancelled by the weights of fixed vertices.

  To prove (\ref{equa:Zht_impair_ax}), we multiply the left-hand side by
  \begin{displaymath}
    \prod_{N+1\leq k\leq 2N} \sigma(a^2 x_k \overline{y}),
  \end{displaymath}
  which is equivalent to adding to the line of parameter $y$ a new line
  $\overline{a}y$ just below the grid;
  Lemma~\ref{lemm:tour_grille} transforms the graph of
  Figure~\ref{fig:Zht_impair_ax}(a) into the graph of Figure \ref{fig:Zht_impair_ax}(b). 
  When we put  $y=ax_1$, we get the indicated fixed edges, which gives as partition function
  \begin{displaymath}
    \prod_{N+1\leq k\leq 2N} \sigma^2(ax_k\overline{x}_1)
    \prod_{1\leq k\leq N} \sigma(a^2x_1\overline{x}_k)
    Z_{\textsc{HT}}(2N;X_{2N}\backslash x_1, x_1,x).
  \end{displaymath}

  \begin{figure}[htbp]
    \begin{center}
      \psset{unit=5mm}
      \begin{pspicture}[.5](-1,-1)(9,11)
	\rput(1,1){\IceGrid{4}{10}}
	\rput(1,0){\linB{4}}
	\rput(1,1){\linB{4}}
	\rput(0,2){\colD{9}}
	\rput(1,10){\linH{4}}
	\psline(5,2)(5,5)
	\multido{\i=1+1}{4}{\rput(\i,1){\flecheG}}
	\multido{\i=2+1}{4}{\psline(4,\i)(5,\i)}
	\multido{\i=7+1}{4}{\psline(4,\i)(5,\i)}
	\multido{\i=1+1}{4}{\psarc(5,6){\i}{270}{90}}
	\psarc(4,2){1}{270}{0}
	\psarc(4,5){1}{0}{90}
	\SpecialCoor
	\rput(4,5){\psline[linestyle=dotted](.5;45)(1.5;45)}
	\rput(4,2){\psline[linestyle=dotted](.5;315)(1.5;315)}
	\rput[r](0,2){$x_1$}
	\rput[r](0,5){$x_{N}$}
	\rput[r](0,6){$x$}
	\rput[r](0,1){$\overline{a}y$}
	\rput[t](1,-.1){$x_{N+1}$}
	\rput[t](4,-.1){$x_{2N}$}
	\rput[l](5.1,1.2){$y$}
	\rput[c](4.5,-1.5){(a)}
      \end{pspicture} \ \ 
      \begin{pspicture}[.5](-1,-1)(10,11)
	\rput(2,1){\IceGrid{4}{10}}
	\rput(1,0){\linB{5}}
	\rput(2,1){\linB{4}}
	\rput(2,9){\linH{4}}
	\rput(2,10){\linH{4}}
	\rput(0,1){\colD{4}}
	\rput(1,2){\colD{3}}
	\rput(1,6){\colD{5}}
	\multido{\n=.5+1.0}{5}{\psarc(5,5.5){\n}{270}{90}}
	\multido{\i=1+1}{3}{\rput(1,\i){\flecheH}}
	\multido{\i=2+1}{4}{\rput(\i,1){\flecheG}}
	\multido{\i=2+1}{3}{\rput(\i,10){\flecheD}}
	\psarc(2,4){1}{90}{180}
	\psarc{<-}(5,5.5){4.5}{0}{90}
	\SpecialCoor
	\rput(2,4){\psline[linestyle=dotted](.5;135)(1.5;135)}
	\rput(5,5.5){\psline[linestyle=dotted](.2;0)(1;2)}
	\rput[r](0,1){$x_1$}
	\rput[r](0,4){$x_N$}
	\rput[r](1,6){$x$}
	\rput[t](1,-.1){$y$}
	\rput[t](2.5,-.1){$x_{N+1}$}
	\rput[t](5,-.1){$x_{2N}$}
	\rput[b](.8,5){$\overline{a}y=x_1$}
	\rput[c](5,-1.5){(b)}
      \end{pspicture}
    \end{center}
    \caption{Proof of (\ref{equa:Zht_impair_ax})}\label{fig:Zht_impair_ax}
  \end{figure}

  Since
  $a^2x_k\overline{y}=ax_k\overline{x}_1$, the equation simplifies. 
  To conclude, we observe that if we start with an edge going out from the crossing $x/x_{2N}$ 
  (function $Z_{\textsc{HT}}^{\downright}$) we get at the end the same orientation (function
  $Z_{\textsc{HT}}^{\downc}$).

The proof of equations (\ref{equa:Zht_impair_bax}-\ref{equa:Zht_pair_bax}) follows the same path. We give in Figure \ref{fig:Zht_impair_bax} the modifications performed to the graph for the proof of (\ref{equa:Zht_impair_bax}).

  \begin{figure}[htbp]
    \begin{center}
      \psset{unit=5mm}
      \begin{pspicture}[.5](-1,-1)(9,11)
	\rput(1,1){\IceGrid{4}{10}}
	\rput(1,0){\linB{4}}
	\rput(0,1){\colD{9}}
	\rput(1,9){\linH{4}}
	\rput(1,10){\linH{4}}
	\psline(5,6)(5,9)
	\multido{\i=1+1}{4}{\rput(\i,10){\flecheG}}
	\multido{\i=1+1}{4}{\psline(4,\i)(5,\i)}
	\multido{\i=6+1}{4}{\psline(4,\i)(5,\i)}
	\multido{\i=1+1}{4}{\psarc(5,5){\i}{270}{90}}
	\psarc(4,6){1}{270}{0}
	\psarc(4,9){1}{0}{90}
	\SpecialCoor
	\rput(4,6){\psline[linestyle=dotted](.5;315)(1.5;315)}
	\rput(4,9){\psline[linestyle=dotted](.5;45)(1.5;45)}
	\rput[r](0,1){$x_1$}
	\rput[r](0,4){$x_N$}
	\rput[r](0,5){$x$}
	\rput[r](0,10){$ay$}
	\rput[t](1,-.1){$x_{N+1}$}
	\rput[t](4,-.1){$x_{2N}$}
	\rput[l](5.1,9.7){$y$}
	\rput[c](4.5,-1){(a)}
      \end{pspicture}\ \ 
      \begin{pspicture}[.5](-1,-1)(10,11)
	\rput(2,1){\IceGrid{4}{10}}
	\rput(2,0){\linB{4}}
	\rput(2,1){\linB{4}}
	\rput(1,1){\colD{5}}
	\rput(0,7){\colD{4}}
	\rput(1,7){\colD{3}}
	\rput(1,10){\linH{5}}
	\rput(2,9){\linH{4}}
	\multido{\i=8+1}{3}{\rput(1,\i){\flecheB}}
	\multido{\i=2+1}{3}{\rput(\i,1){\flecheD}}
	\multido{\i=2+1}{4}{\rput(\i,10){\flecheG}}
	\psarc(2,7){1}{180}{270}
	\multido{\n=.5+1.0}{5}{\psarc(5,5.5){\n}{270}{90}}
	\psarc{->}(5,5.5){4.5}{270}{0}
	\SpecialCoor
	\rput(2,7){\psline[linestyle=dotted](.5;225)(1.5;225)}
	\rput(5,5.5){\psline[linestyle=dotted](.2;0)(1.2;0)}
	\rput[r](1,1){$x_1$}
	\rput[r](1,4){$x_N$}
	\rput[r](1,5){$x$}
	\rput[t](2,-.1){$x_{N+1}$}
	\rput[t](5,-.1){$x_{2N}$}
	\rput[b](1,11){$y=\overline{a}x_{1}$}
	\rput[r](1.5,5.5){$ay=x_1$}
	\rput[c](5,-1){(b)}
      \end{pspicture}
    \end{center}
    \caption{Proof of (\ref{equa:Zht_impair_bax})}
    \label{fig:Zht_impair_bax}
  \end{figure}

\endproof

\begin{lemm}{\em [specialization of $Z_{\textsc{QT}}$]}
\label{lemm:specZQT}
If we denote
\begin{align*}
  \overline{A}_{Q}(x_1,X_{m-1}\backslash x_1) & =  
  \prod_{1\leq k\leq m-1} \sigma(a^2 x_k\overline{x}_1) 
   \sigma(ax_1 \overline{x}_k),\\
  A_{Q}(x_1;X_{m-1}\backslash x_1) & = 
  \prod_{1\leq k\leq m-1} \sigma(a^2 x_1 \overline{x}_k)
\sigma(a x_k \overline{x}_1),
\end{align*}
then for $\star = \convcorner,\divcorner,\downright,\upleft$ and
$\square=\upleft,\downright,\convcorner,\divcorner$ respectively, we have:
\begin{align}
\!\!\!\!\!\!  Z_{\textsc{QT}}^{\star}(2m;X_{m-1},\mathbf{\overline{a}x_1},y) &=
  \sigma(ax_1\overline{y}) \overline{A}_{Q}(x_1,X_{m-1})
  Z_{\textsc{QT}}^{\square}(2m-2;X_{m-1}\backslash x_1, y,x_1),
  \label{eq:Zqt_spe_bax}\\
\!\!\!\!\!\!   Z_{\textsc{QT}}^{\square}(2m;X_{m-1},x,\mathbf{ax_1}) & = 
\sigma(ax\overline{x}_1) A_{Q}(x_1;X_{m-1}\backslash x_1)
  Z_{\textsc{QT}}^{\star}(2m-2;X_{m-1}\backslash x_1, x_1,x).
  \label{eq:Zqt_spe_ax}
\end{align}
\end{lemm}

\proof
The proof is very similar to the previous one: we add weighted vertices
to the graph, then apply Lemma~\ref{lemm:tour_grille}, identify the edges 
that are fixed by the specialization of the parameter, and conclude
by observing that the identity that we obtain can be simplified by
the added weights. The corresponding graphs are given in
Figure~\ref{fig:Zqt_spe}; the symbol $\triangle$ means a change of orientation
only when $m$ is even.

  \begin{figure}[htbp]
    \begin{align*}
    \psset{unit=.5cm}
    \begin{pspicture}[.5](-1,-1)(12,11)
      \rput(2,1){\IceGrid{5}{6}}
      \rput(2,0){\linB{6}}
      \rput(1,1){\colD{5}}
      \rput(0,1){\colD{5}}
      \multido{\i=0+1}{5}{\rput(1,\i){\flecheH}}
      \multido{\i=1+1}{5}{\psline(6,\i)(7,\i)}
      \multido{\i=1+1}{5}{\psarc(7,6){\i}{270}{180}}
      \psline(7,1)(7,5)
      \psarc(6,5){1}{0}{90}
      \psarc(2,5){1}{90}{180}
      \SpecialCoor
      \rput(7,6){\psdots[dotsize=1mm](1;45)(2;45)(3;45)(4;45)(5;45)}
      \rput(6,5){\psline[linestyle=dotted](.5;45)(1.5;45)}
      \rput(2,5){\psline[linestyle=dotted](.5;135)(1.5;135)}
      \rput(6,5){\psdot[dotstyle=triangle,dotsize=1.5mm](1;45)}
      \rput[r](0,1){$x_1$}
      \rput[r](0,5){$x_{m-1}$}
      \rput[b](1.2,6){$x$}
      \rput[t](1,-.1){$ax$}
      \rput[t](7,-.1){$y$}
    \end{pspicture}
    & =  
    \psset{unit=.5cm}
    \begin{pspicture}[.5](-1,-1)(12,11)
      \rput(1,2){\IceGrid{7}{5}}
      \rput(0,1){\colD{6}}
      \rput(1,0){\linB{5}}
      \rput(2,1){\linB{4}}
      \multido{\i=1+1}{5}{\rput(1,\i){\flecheH}}
      \rput(1,2){\colD{5}}
      \psarc{<-}(7,6.5){4.5}{60}{90}
      \psarc{->}(7,6.5){4.5}{0}{30}
      \rput(2,2){\linB{6}}
      \multido{\i=2+1}{5}{\rput(\i,2){\flecheD}}
      \multido{\i=1+1}{5}{\psarc(6,6){\i}{90}{180}}
      \multido{\n=.5+1.0}{5}{\psarc(7,6.5){\n}{270}{90}}
      \multido{\i=7+1}{5}{\psline(6,\i)(7,\i)}
      \multido{\i=2+1}{4}{\rput(\i,1){\flecheG}}
      \psarc(5,2){1}{270}{0}
      \rput(7,2){\flecheB}
      \psarc(6.5,6){.5}{0}{180}
      \SpecialCoor
      \rput(7,6.5){\psdots[dotsize=1mm](.5;45)(1.5;45)(2.5;45)(3.5;45)(4.5;45)}
      \rput(6.5,6){\psline[linestyle=dotted](.2;90)(.8;90)}
      \rput(5,2){\psline[linestyle=dotted](.5;315)(1.5;315)}
      \rput(6.5,6){\psdot[dotstyle=triangle,dotsize=2mm](.5;90)}
      \rput[r](0,1){$x=\overline{a}x_1$}
      \rput[t](1,-.1){$x_1$}
      \rput[r](5,-.1){$x_{m-1}$}
      \rput[t](7,.9){$y$}
      \rput[l](6,1){$ax$}
    \end{pspicture}\\
    \psset{unit=.5cm}
    \begin{pspicture}[.5](-1,-1)(11,12)
      \rput(1,1){\IceGrid{5}{7}}
      \rput(1,0){\linB{5}}
      \rput(0,2){\colD{6}}
      \rput(1,1){\linB{5}}
      \multido{\i=1+1}{5}{\rput(\i,1){\flecheG}}
      \multido{\i=2+1}{5}{\psline(5,\i)(6,\i)}
      \psline(6,2)(6,6)
      \psarc(5,2){1}{270}{0}
      \psarc(5,6){1}{0}{90}
      \multido{\i=1+1}{5}{\psarc(6,7){\i}{270}{180}}
      \SpecialCoor
      \rput(5,2){\psline[linestyle=dotted](.5;315)(1.5;315)}
      \rput(5,6){\psline[linestyle=dotted](.5;45)(1.5;45)}
      \rput(6,7){\psdots[dotsize=1mm](1;45)(2;45)(3;45)(4;45)(5;45)}
      \rput(5,6){\psdot[dotstyle=triangle,dotsize=2mm](1;45)}
      \rput[r](0,1){$\overline{a}y$}
      \rput[r](0,2){$x_1$}
      \rput[r](0,6){$x_{m-1}$}
      \rput[r](0,7){$x$}
      \rput[l](6.2,1){$y$}
    \end{pspicture}
    & =  
    \psset{unit=5mm}
    \begin{pspicture}[.5](-1,-1)(11,12)
      \rput(2,1){\IceGrid{5}{7}}
      \rput(1,0){\linB{6}}
      \rput(0,1){\colD{5}}
      \rput(1,2){\colD{4}}
      \rput(1,7){\flecheD}
      \multido{\i=1+1}{4}{\rput(1,\i){\flecheH}}
      \multido{\i=2+1}{5}{\rput(\i,1){\flecheG}}
      \rput(2,2){\colD{6}}
      \rput(2,1){\linB{5}}
      \multido{\i=2+1}{6}{\rput(2,\i){\flecheB}}
      \multido{\i=1+1}{5}{\psarc(6,6){\i}{270}{0}}
      \multido{\i=7+1}{5}{\psline(\i,6)(\i,7)}
      \multido{\n=.5+1.0}{5}{\psarc(6.5,7){\n}{0}{180}}
      \psarc(2,5){1}{90}{180}
      \psarc(6,6.5){.5}{270}{90}
      \SpecialCoor
      \rput(6.5,7){\psdots[dotsize=1mm](.5;45)(1.5;45)(2.5;45)(3.5;45)(4.5;45)}
      \rput(6,6.5){\psline[linestyle=dotted](.2;0)(.8;0)}
      \rput(6,6.5){\psdot[dotstyle=triangle,dotsize=2mm](.5;0)}
      \rput(2,5){\psline[linestyle=dotted](.5;135)(1.5;135)}
      \psarc{<->}(6.5,7){4.5}{30}{60}
      \rput[r](1,7){$x$}
      \rput[r](0,1){$x_1$}
      \rput[r](0,5){$x_{m-1}$}
      \rput[t](1,-.1){$y=ax_1$}
      \rput[b](1.1,6){$\overline{a}y$}
    \end{pspicture}
    \end{align*}
    \caption{Proof of (\ref{eq:Zqt_spe_bax}-\ref{eq:Zqt_spe_ax})}
    \label{fig:Zqt_spe}
  \end{figure}
\endproof

\begin{rema}\label{rema:sym}
By using the (pseudo-)symmetry in $(x,y)$, we may transform any specialization of the variable $y$ into a specialization of the variable $x$. Moreover, by using Lemma \ref{lemm:sym} and (when $a=\omega_6$) Lemma \ref{lemm:symZomega6}, we obtain for $Z$, $Z_{\textsc{HT}}$  and $Z_{\textsc{QT}}$, $2N$ independent specializations of the variable $x$. 
\end{rema}

\subsection{Special value of the parameter $a$; conclusion}

When $a=\omega_6=\exp(i\pi/3)$, two new ingredients may be used. The first one is Lemma \ref{lemm:symZomega6}, as mentioned in Remark \ref{rema:sym}. The second one is that with this special value of $a$ we have:
\begin{equation}\label{spec}
\si(a)=\si(a^2) \ \ \ \ \ \ \ \si(a^2x)=-\si(\ba x)=\si(a\bx).
\end{equation}
which implies that the products appearing in Lemmas \ref{lemm:specZ}, \ref{lemm:specZHT} and \ref{lemm:specZQT} may be written in a more compact way:
\begin{align*}
  A(x_{N+1},X_{2N}\backslash\{x_1,x_{N+1}\}) & =  \sigma(a)
  \prod_{k\neq 1,N+1} \sigma(a x_k \overline{x}_{N+1}),\\
  \overline{A}(x_{N+1},X_{2N}\backslash \{x_1,x_{N+1}\}) & =  
  \sigma(a) \prod_{k\neq 1,N+1} \sigma(a x_{N+1} \overline{x}_k),\\
  A_{H}^{1}(x_1,X_{2N}\backslash x_1) & =  
  \prod_{1\leq k\leq 2N} \sigma(a x_{k} \overline{x}_1),\\
  \overline{A}_{H}^{1}(x_1,X_{2N}\backslash x_1) & =  
  \prod_{1\leq k\leq 2N} \sigma(a x_1 \overline{x}_k),\\
  A_{H}^{0}(x_N,X_{2N-1}\backslash x_N) & = 
  \prod_{1\leq k\leq 2N-1} \sigma(a x_{k} \overline{x}_{N}),\\
  \overline{A}_{H}^{0}(x_N,X_{2N-1}\backslash x_{N}) &=
  \prod_{1\leq k\leq 2N-1} \sigma(a x_{N} \overline{x}_k),\\
  \overline{A}_{Q}(x_1,X_{m-1}\backslash x_1) & = 
  \prod_{1\leq k\leq m-1} \sigma^2(a x_1 \overline{x}_k),\\
  A_{Q}(x_1,X_{m-1}\backslash x_1) & =  
  \prod_{1\leq k\leq m-1} \sigma^2(a x_k \overline{x}_1).
\end{align*}

Thus we get by comparing:
\begin{align*}
  A(x_i, X_{2N}\backslash x_i, x)
  A_{H}^{1}(x_i,X_{2N}\backslash x_i)& = \sigma(a x
  \overline{x}_i) A_{Q}(x_i,X_{2N}\backslash x_i)\\
  \overline{A}(x_i,X_{2N}\backslash x_i,x)
  \overline{A}_{H}^{1}(x_i,X_{2N}\backslash x_i) & = 
  \sigma(a x_i \overline{x})
  \overline{A}_{Q}(x_i,X_{2N}\backslash x_i),
\end{align*}
whence (\ref{eq:main1.1}) and (\ref{eq:main1.2}) imply that
(\ref{eq:main2.1}) and (\ref{eq:main2.2}) are true (in size
$4N+2$) for the $2N$ specializations $x=a^{\pm 1}x_i$ ($1\leq i\leq
N$). It is enough to prove (\ref{eq:main2.2}) (Laurent polynomials of half-width
 $2N-1$), but we still need one specialization to get
 (\ref{eq:main2.1}) (half-width $2N$).

For (\ref{eq:main1.1}) and (\ref{eq:main1.2}), we observe the same kind of simplification 
\begin{displaymath}
  A(x_i, X_{2N-1}\backslash x_i) \sigma(a x \overline{x}_i)
  A_{H}^{0}(x_i,X_{2N-1}\backslash x_i) = \sigma(a x \overline{x}_i)
  A_{Q}(x_i,X_{2N-1}\backslash x_i),
\end{displaymath}
whence (\ref{eq:main2.2}) and (\ref{eq:main2.1}) for the size $4N-2$
imply that (\ref{eq:main1.1}) and (\ref{eq:main1.2}) are true for the 
 $N$ specializations $x=ax_i$, $N\leq i\leq 2N-1$. We obtain in the same way
the coincidence for the $N$ specializations $x=\overline{a}x_i$, $N\leq i\leq
2N-1$. Thus we have $2N$ specialiations of $x$: it is enough both for
 (\ref{eq:main1.1}) (half-width  $2N-1$), and for (\ref{eq:main1.2}) (half-width $2N-2$).

At this point, we have \emph{almost} proved
\begin{center}
  ((\ref{eq:main1.1}) and (\ref{eq:main1.2}), in size $4N$) $\Longrightarrow$
  ((\ref{eq:main2.1}) and (\ref{eq:main2.2}), in size $4N+2$) $\Longrightarrow$
  ((\ref{eq:main1.1}) and (\ref{eq:main1.2}), in size $4N+4$);
\end{center}
\emph{almost}, because we still need {\em one} specialization for
(\ref{eq:main2.1}).

We get this missing specialization, not directly for 
$\ZqDownright$, $\ZqUpleft$, $\ZhDownright$ and $\ZhUpleft$, but for the original series $Z_{\textsc{QT}}(4N+2;X_{2N},x,y)$ and
$Z_{\textsc{HT}}(2N+1;X_{2N},x,y)$: indeed if we set
$x=ay$ we may apply
Lemma~\ref{lemm:tour_grille}.

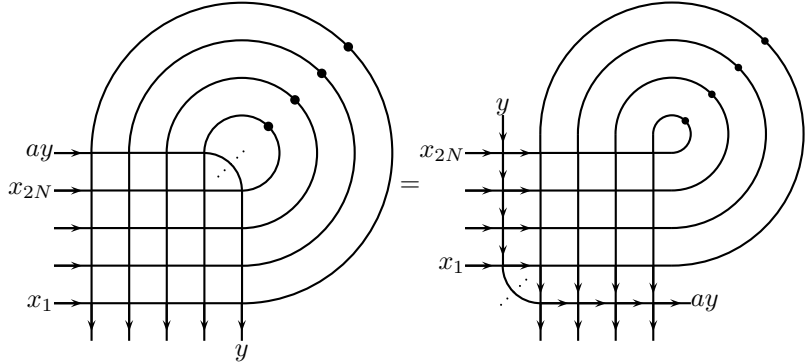
$$  \psset{unit=5mm}
  \begin{pspicture}[.5](-1,-1)(9,9)
    \rput(0,0){\qQTIce{5}}
    \rput[r](0,1){$x_1$}
    \rput[r](0,4){$x_{2N}$}
    \rput[r](0,5){$ay$}
    \rput[t](5,-.1){$y$}
  \end{pspicture}
 =  
  \psset{unit=5mm}
  \begin{pspicture}[.5](-1,-1)(9,9)
    \rput(2,2){\IceGrid{4}{4}}
    \rput(2,1){\linB{4}}
    \rput(2,0){\linB{4}}
    \rput(1,2){\colD{4}}
    \rput(0,2){\colD{4}}
    \multido{\i=2+1}{4}{\rput(\i,1){\flecheD}}
    \multido{\i=3+1}{4}{\rput(1,\i){\flecheB}}
    \multido{\i=2+1}{4}{\psline(\i,5)(\i,5.5)\psline(5,\i)(5.5,\i)}
    \multido{\n=.5+1.0}{4}{\psarc(5.5,5.5){\n}{270}{180}}
    \SpecialCoor
    \rput(2,2){\psline[linestyle=dotted](.5;225)(1.5;225)}
    \rput(5.5,5.5){\psdots[dotsize=1mm](.5;45)(1.5;45)(2.5;45)(3.5;45)}
    \psarc(2,2){1}{180}{270}
    \rput[b](1,6){$y$}
    \rput[r](0,2){$x_1$}
    \rput[r](0,5){$x_{2N}$}
    \rput[l](6,1){$ay$}
  \end{pspicture}$$
$$Z_{\textsc{QT}}(4N+2;X_{2N},{\bf ay},y)  =  
\sigma(a) \prod_{1\leq k\leq 2N} \sigma(a x_k \overline{y}) \sigma(a^2 y \overline{x}_k) Z_{\textsc{QT}}(4N;X_{2N}\backslash x_{2N}, x_{2N},x_{2N})$$
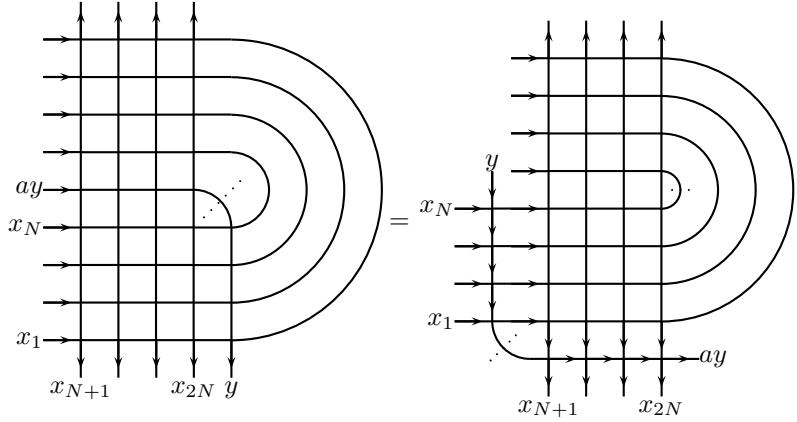
$$  \psset{unit=5mm}
  \begin{pspicture}[.5](-1,-1)(9,9)
    \rput(0,0){\HTIceOdd{4}}
    \rput[r](0,1){$x_1$}
    \rput[r](0,4){$x_N$}
    \rput[r](0,5){$ ay$}
    \rput[t](1,-.1){$x_{N+1}$}
    \rput[t](4,-.1){$x_{2N}$}
    \rput[t](5,-.1){${y}$}
  \end{pspicture}
  =  
  \psset{unit=5mm}
  \begin{pspicture}[.5](-.5,-1)(9,10)
    \rput(2,1){\HTIceEven{4}}
    \rput(3,0){\linB{4}}
    \rput(.5,2){\colD{4}}
    \multido{\i=2+1}{4}{\psline(1.5,\i)(2,\i)}
    \psline(2.5,1)(3,1)
    \multido{\i=3+1}{4}{\rput(\i,1){\flecheD}}
    \multido{\i=3+1}{4}{\rput(1.5,\i){\flecheB}}
    \SpecialCoor
    \psarc(2.5,2){1}{180}{270}
    \rput(2.5,2){\psline[linestyle=dotted](.5;225)(1.5;225)}
    \rput[r](.5,2){$x_1$}
    \rput[r](.5,5){$x_{N}$}
    \rput[b](1.5,6){$y$}
    \rput[t](3,-.1){$x_{N+1}$}
    \rput[t](6,-.1){$x_{2N}$}
    \rput[l](7,1){$ay$}
  \end{pspicture}$$
$$  Z_{\textsc{HT}}(2N+1;X_{2N},{\bf ay},y)  = 
\left(\prod_{1\leq k\leq N} \sigma(a x_k \overline{y}) 
\prod_{N+1\leq k\leq 2N} \sigma(a^2 y \overline{x}_k)\right)
  Z_{\textsc{HT}}(2N; X_{2N}\backslash x_{N}, x_N,x_N)
$$

This way, we get another point where
(\ref{eq:main2}) is true, and thus, because we already have
(\ref{eq:main2.2}), by difference we obtain that (\ref{eq:main2.1}) 
holds for $y=\overline{a}x$.

This completes the proof of Theorem \ref{theo:main}.

\section{Open questions}\label{sec:open}

The first open problem concerns the so-called $q$-enumeration of
ASM's, which consists in counting ASM's (or classes of symmetry of
ASM's) with respect to their number of $-1$ entries (to the number of
orbits of $-1$ entries in the case of symmetric ASM's).  For a generic
value of the global parameter $a$, when we put all variables to $1$,
the weight of zero entries in the ASM ($\si(a)$) is different from the
weight of non-zero entries ($\si(a^2)$). This may allow a
$q$-enumeration of ASM's.  But in our case, since we fix the value of
$a$ to $\exp(i\pi/3)$, we have $\si(a)=\si(a^2)$, thus we cannot keep
the trace of non-zero entries. The partition functions as we have
defined them do not seem to factor for other values of $a$.

The second question is a very frustrating one, and is common to all
these beautiful equinumeration formulas: is it possible to give a
bijective explanation to equation (\ref{eq:form}), as well as to
equations (\ref{eq:kup}) and (\ref{eq:QT_m1}-\ref{eq:QT_p1})? 

As an indication of what the ``right'' bijection should look like, let
us note that, in the qQTASM case, the restricted ice models that
define the partition functions $\ZqUpleft(4N+2)$ and
$\ZqDownright(4N+2)$ correspond respectively to those qQTASMs where
the four central entries are $(0,-1,-1,0)$ and $(1,0,0,1)$; similarly,
the ice models that defins partition functions $\ZhUpleft(2N+1)$ and
$\ZhDownright(2N+1)$ correspond respectively to those HTASMs where the
center entry is $1$ and $-1$. Thus, a consequence of Equations
\pref{eq:main2.1} and \pref{eq:main2.2} is that the total proportion
of qQTASMs (of size $4N+2$) with two negative entries among the four
central entries, is exactly the proportion of HTASMs (of size $2N+1$)
with negative central entry; in \cite{RSHT}, Razumov and Stroganov proved
this proportion to be exactly $\frac{N}{2N+1}$. This observation gives a new
occurrence of the {\em $1/N$ phenomenon}, as defined in \cite{stro}.

A last question deals with the {\em link pattern} distribution of
FPL's (Fully Packed Loop configurations, in bijection with ASM's).
This question arose in the intriguing Razumov-Stroganov conjecture
\cite{RSConj1,RSConj2}. It appears that Equations (\ref{eq:form}) and
(\ref{eq:kup}) can be refined into
Conjectures 5 and 6 presented in \cite{duchon}.



\begin{thebibliography}{10}

\bibitem{aval}
{\sc J.-C. Aval},
{\em On the symmetry of the partition function of some square ice models},
to appear in Theor. Math. Phys, {\tt arXiv:0903.0777}.

\bibitem{duchonHDR}
{\sc P. Duchon},
{\em Configurations de boucles compactes},
Habilitation \`a diriger des recherches, Universit\'e Bordeaux 1, 2009 (in French).

\bibitem{duchon}
{\sc P. Duchon},
{\em On the link pattern distribution of quarter-turn symmetric FPL configurations},
Proceedings of FPSAC 2008, 12p.

\bibitem{izergin}
{\sc A. G. Izergin},
{\em Partition function of the six-vertex model in a finite volume},
Sov. Phys. Dokl. {\bf 32} (1987) 878--879.

\bibitem{korepin}
{\sc V. E. Korepin},
{\em Calculation of Norms of Bethe Wave Functions},
Comm. Math. Phys. {\bf 86} (1982) 391--418.

\bibitem{kup1}
{\sc G. Kuperberg},
{\em Another proof of the alternating sign matrices conjecture},
Internat. Math. Research Not. ,{\bf 1996} (1996), 139--150. 

\bibitem{kup}
{\sc G. Kuperberg},
{\em Symmetry classes of alternating sign matrices under one roof},
Ann. Math. {\bf 156} (2002) 835--866.

\bibitem{MRR}
{\sc W. Mills, D. Robbins, H. Rumsey},
{\em Alternating sign matrices and descending plane partitions}, J. Combin. Th. Ser. A {\bf 34} (1983), 340--359.

\bibitem{RSConj1}
{\sc A. V. Razumov, Y. G. Stroganov},
{\em Spin chains and combinatorics},
J. Phys. A, {\bf 34} (2001), 3185--3190.

\bibitem{RSConj2}
{\sc A. V. Razumov, Y. G. Stroganov},
{\em Combinatorial nature of the ground-state vector of the o(1) loop model} 
Theoret. and Math. Phys., {\bf 138} (2004),333--337.

\bibitem{RSHT}
{\sc A. V. Razumov, Y. G. Stroganov},
{\em Enumeration of half-turn symmetric alternating sign matrices of odd order},
Theor. Math. Phys. {\bf 148} (2006) 1174--1198.

\bibitem{RS}
{\sc A. V. Razumov, Y. G. Stroganov},
{\em Enumeration of quarter-turn symmetric alternating sign matrices of odd order},
Theoret. Math. Phys. {\bf 149} (2006) 1639--1650.

\bibitem{robbins}
{\sc D. Robbins},
{\em Symmetry classes of alternating sign matrices}, {\tt arXiv:math.CO/0008045}.

\bibitem{IKdet}
{\sc Y. G. Stroganov},
{\em A new way to deal with Izergin-Korepin determinant at root of unity},
{\tt arXiv:math-ph/0204042}.

\bibitem{stro}
{\sc Y. G. Stroganov},
{\em $1/N$ phenomenon for some symmetry classes of the odd alternating sign matrices},
{\tt arXiv:0807.2520}.


\bibitem{zeil}
{\sc D. Zeilberger},
{\em Proof of the alternating sign matrix conjecture}, Electronic J. Combinatorics {\bf 3}, No. 2 (1996)
, R13, 1--84.


\end{thebibliography}
\end{document}